\input amstex
\documentstyle{amsppt}
\mag=\magstep1
\pagewidth{13cm}
\NoBlackBoxes
\pageheight{18cm}
\def\a{\alpha}
\def\ol{\overline}
\def\b{\beta}

\def\gam{\gamma}
\def\Gam{\Gamma}

\def\lam{\lambda}
\def\ome{\omega}
\def\Ome{\Omega}
\def\sig{\sigma}

\def\A{\Cal A }
\def\B{\Cal B}
\def\E{\Cal E}
\def\bE{\bold E}

\def\G{G_2}

\def\Diff{\text{\rm Diff}}
\def\R{\Bbb R}
\def\C{\Bbb C}
\def\H{\text{\rm H}}
\def\Z{\Bbb Z}
\def\O{\Bbb O}
\def\P{\frak P}
\def\L{\Cal L}
\def\M{\frak M}

\def\w{\wedge}
\def\({\left(}
\def\){\right)}
\def\G{G_2}

\def\CL{{\text{ \rm CL}}}

\def\h{\hat}
\def\hrho{\hat{\rho}}

\def\til{\tilde}
\def\wtil{\widetilde}

\def\l{\left}
\def\pa{\partial}
\def\olpa{\ol{\partial}}
\def\r{\right}
\def\ran{\rangle}
\def\lan{\langle}
\def\ss{\scriptscriptstyle}
\def\trian{\triangle}

\def\GL{{\text{\rm GL}}}
\def\SL{{\text{\rm SL}}}
\def\Spin{\text{\rm Spin}}
\def\Pin{\text{\rm Pin}}
\def\CSpin{\text{\rm CSpin}}
\def\CPin{\text{\rm Cpin}}

\def\SU{\text{\rm SU}}
\def\SO{\text{\rm SO}}

\def\I{\Cal I}
\def\J{\Cal J}
\def\K{\Cal K}
\def\VV{V\oplus V^*} 
\def\TT{T\oplus T^*}
\def\CY{{\text{\rm CY}}}
\def\HK{\text{\rm HK}}
\def\CL{\text{\rm CL}}
\def\even{\text{\rm even}} 
\def\odd{\text{\rm odd}}
\def\Ad{\text{\rm Ad}}
\def\t{\theta}
\def\Ob{\text{\rm Ob}}
\def\vol{\text{\rm vol}}
\document

\topmatter
\title\nofrills\bf
On  deformations \\of \\generalized Calabi-Yau, hyperK\"ahler, 
$\G$ and Spin$(7)$ structures I
\endtitle
\rightheadtext{On deformations of generalized structures}
\leftheadtext{ Ryushi Goto}
\author{  Ryushi Goto}\endauthor
\affil {Department of Mathematics, \\Graduate School of Science,\\ Osaka university
 }\endaffil
 \address 1-1Toyonaka, Osaka, 560 Japan 
 \endaddress
\email goto{\@\,}math.sci.osaka-u.ac.jp \endemail
\abstract
In this paper we will introduce a new notion of  geometric structures defined by systems of closed differential forms in term of the Clifford algebra of the direct sum of the tangent bundle and the cotangent bundle on a manifold. 
We develop a unified approach of a deformation problem and establish a criterion of 
unobstructed deformations of the structures from a cohomological point of view. 
We construct the moduli spaces of the structures by using the action of $b$-fields and show that the period map of the moduli space is locally injective under a cohomological condition
(the local Torelli type theorem).
We apply our approach to generalized Calabi-Yau (metrical) structures and obtain an analog of the theorem by Bogomolov-Tian-Todorov.
Further we prove that deformations of generalized Spin$(7)$ structures are unobstructed.
\endabstract

\toc
\head \S 0 Introduction\endhead
\head \S 1 Clifford algebras and Spin groups\endhead
\head \S 2 Clifford-Lie operators \endhead
\head \S 3 Deformations of generalized geometric structures \endhead 
\head \S 4 Deformations of generalized Calabi-Yau structures \endhead 
\head \S5. Generalized hyperK\"ahler, $\G$ and Spin$(7)$ structures 
\endhead
\endtoc
\endtopmatter

\head \S0 Introduction \endhead
In a paper [Go], the author introduced a notion of geometric structures defined by 
systems of closed differential forms which are based on the action of the gauge group of the tangent bundle on a manifold. 
This approach provides a systematic construction of smooth moduli spaces of Calabi-Yau, hyperK\"ahler, G$_2$ and Spin$(7)$
structures. In a paper [Hi1] Hitchin presented a generalized geometry, which depend on
a suggestive idea of replacing the tangent bundle by the direct sum of the tangent bundle $T$ and 
the cotangent bundle $T^*$.
The generalized geometry is a current issue which is rapidly 
studied in differential geometry and mathematical physics. 
However the idea of generalized geometry is perhaps not as widely appreciated as it should be. 
A possible reason for this is that the generalized geometry is at present
restricted to rather special cases.
In this paper we will develop the idea of the generalized geometry from a wide view point as in [Go] which is of general nature together with some new applications.
Since there is an indefinite metric on the direct sum $\TT$ on a manifold $X$, then the bundle of the Clifford algebras $\CL(X)$ of $\TT$ naturally appears  and
we obtain various fibre bundles of Lie groups such as Spin, Pin and conformal Pin group $\CPin(X)$ which act on the differential forms on $X$ by the spin representation.
Then we introduce a notion of geometric structures which are defined by systems of closed differential forms in an orbit of the action of the conformal pin group. 
We develop a deformation problem of the structures and establish a criterion for unobstructed deformations of the structures and study a problem for when the local Torelli type theorem holds (theorem 3-7, 8 and 9). 
Then we apply our approach to interesting special cases of generalized structures 
discussed in [Hi1], [Gu1].
For instance, generalized SL$_n(\C)$ structures are defined as complex pure spinors with 
a non-degenerate condition which are a generalization of complex structures with trivial canonical line bundle. 
Note that we call them generalized $SL_n(\C)$ structures because the special linear group SL$_n(\C)$ naturally arises as the isotropy group. (see [Go] for SL$_n(\C)$ structures).
A generalized SL$_n(\C)$ structure $\phi$ induces a generalized complex structure $\J_\phi$.  Then our criterion easily implies 
\proclaim{Theorem 4-1-6} 
If we have the $dd^\J$-property for the $\J_\phi$ corresponding a generalized $\SL_n(\C)$ structure $\phi$, then 
$\phi$ is a topological structure, so that is, we have unobstructed deformations of $\phi$
on which the local Torelli type theorem holds.
\endproclaim
The $dd^\J$- property is a generalization of ordinary $\pa\ol\pa$-lemma in K\"ahler geometry and
Gualtieri shows that the $dd^\J$ property holds for generalized K\"ahler structures [Gu2]. 
A Calabi-Yau (metrical) structure in [Gu1]  is a pair consisting of generalized SL$_n(\C)$ structures $\phi_0$ and $\phi_1$ 
such that
the corresponding pair of generalized complex structures yields
a generalized K\"ahler structure. Deformations of such pairs seems to be complicated, however it is observed that our systematic approach is adapted to obtain unobstructed deformations of Calabi-Yau (metrical) structures and the local Torelli type theorem of them (theorem 4-2-3). 

Li also shows a result of deformations of generalized complex structure [Li].
It is worthwhile to mention that there is a relation between deformations of generalized SL$_n(\C)$ structures and ones of 
generalized complex structures (Proposition 4-1-7). 
If we have the $dd^\J$-property, there is a surjective map from deformations of generalized SL$_n(\C)$ structures to ones of 
generalized complex structures, so that is, 
both deformations are essentially same,
which yields an another proof of the result by Li. 
 We give a brief outline of this paper.
In section 0, we present an exposition of the Clifford algebras of the direct sum of a real vector space $V$ and the dual space $V^*$ 
and various groups such as Spin, Pin and conformal pin groups. 
It is important that the exponential $e^b$ (resp. $e^\b$) for a $2$-form
$b\in\w^2 V^*$ (resp. a $2$-vector $\b\in \w^2V)$ gives  an element of the spin
group.  The materials in this section are 
 already well explained in [L-M], [Ha] and [Hi1].
In section 2, we introduce a subbundle $\CL^k$ over a manifold $X$ which gives a filtration of the even Clifford bundle and one of the odd Clifford bundle : 
$$\align 
&\CL^0\subset \CL^2\subset \CL^4\subset\cdots,\\
&\CL^1\subset\CL^3\subset\CL^5\subset\cdots.
\endalign
$$
Further we discuss 
differential operators acting on differential forms on $X$  which arise as commutators between the exterior derivative $d$ and the action of the Clifford algebra $\CL$. 
The Clifford-Lie operators of order $3$ are introduced in definition 2-2 which are
invariant under the adjoint action of $e^a$ for $a\in \CL^2$ (lemma 2-4). The
exterior derivative $d$ is a Clifford-Lie operator of order $3$ and  it follows
that  
$e^{-a}\circ d\circ e^a$
is also a Clifford-Lie operator of order $3$, which play a significant role in studying the deformation problem.
In section 3, a notion of geometric structures is introduced. 
We start with the direct sum of the real vector space $V$ of $n$ dim and the dual space $V^*$. The conformal pin group $\CPin(\VV)$ of $\VV$ linearly acts on the direct sum of skew-symmetric tensors $\oplus^l\w^*V^*$.   
Let $\Phi=(\phi_1, \cdots, \phi_l)$ be an element of the direct sum $\oplus^l\w^*V^*$ and $\B(V)$ the orbit of $\CPin(\VV)$ through $\Phi$. 
We fix the orbit $\B(V)$ and goes to a oriented, compact manifold $X$ of dim $n$.
The orbit $\B(V)$ yields the orbit in $\oplus^l\w^*T^*_xX$ for  each point $x\in X$ and we have a fibre bunele $\B(X)$ by
$$
\B(X):=\bigcup_{x\in X}\B(T_xX)\to X.
$$
The set of global sections of $\B(X)$ is denoted by $\E_\B(X)$ and 
then we define a $\B(V)$-structure on $X$ by a $d$-closed section of $\E_\B(X)$.
We denote by $\wtil\M_\B(X)$ the set of $\B(V)$-structures on $X$ :
$$
\wtil\M_\B(X)=\{\, \Phi\in \E_\B(X)\, |\, d\Phi=0\, \}.
$$ 
Then we define a moduli space of $\B(V)$-structures on $X$ by the quotient space :
$$
\M_\B(X)=\wtil\M_\B(X)/\wtil\Diff_0(X),
$$
where $\wtil\Diff_0(X)$ is an extension of the diffeomorphisms of $X$ by the action of $d$-exact $b$-fields (see definition 3-2).
Since the de Rham cohomology class of $\Phi$ is invariant under the action of $\wtil\Diff_0(X)$, 
we have the Period map :
$$
P_\b\:\M_\B(X)\to H^*_{dR}(X).
$$
In order to discuss deformations of a $\B(V)$ structure $\Phi$, we introduce a suitable deformation complex $\#_\B$ (proposition 3-3) : 
$$\minCDarrowwidth{0.2cm} 
\CD 
0@>>>\bold E^{-1}(X)@>d_{-1}>>\bold E^0(X)@>d_0>>\bold E^1(X)@>d_1>>
\bold E^2(X)@>d_2>>\cdots,
\endCD
$$
Each vector bundle $\bE^{k-1}(X)$ is defined by the action of the Clifford subbundle $\CL^k$ 
of $\Phi$, so that is, $\bE^{k-1}(X)=\CL^k\cdot\Phi$ and the differential operator $d_k$ is the restriction of $d$ to the bundle $\bE^k(X)$. 
An orbit $\B(V)$ is {\it an elliptic orbit} if the deformation complex $\#_\B$ is an elliptic complex.
It is observed that the complex $\#_\B$ is a subcomplex of the direct sum of the de Rham complex and we have the map $p^k_\B$ from the cohomology groups $H^k(\#_\B)$  of the complex 
$\#_\B$ to the direct sum of the de Rham cohomology groups. 
We say a $\B(V)$-structure $\Phi$ is {\it a topological structure} if  the map $p^k_\B$ is injective for $k=1,2$ (definition 3-5).
Our criterion for unobstructed deformations and the local Torelli type theorem is shown in theorem 3-7 :
\proclaim{Theorem 3-7} 
 Let $\B(V)$ be an elliptic orbit and
 $\Phi$ a $\B(V)$-structure on a compact and oriented $n$-manifold $X$.
 If $\Phi$ is a topological structure, then deformations of $\Phi$ are unobstructed and the  deformation space of $\Phi$ is locally embedded into the de Rham cohomology group $H^*_{dR}(X)$.
 In particular, if an orbit $\B(V)$ is elliptic and topological, the period map $P_\B$ of the moduli space 
 $\M_\B(X)$ of $\B(V)$ structures on $X$ is locally injective.
 \endproclaim
In section $4$ we apply our approach to generalized SL$_n(\C)$ structures and generalized Calabi-Yau (metrical) structures. 
In section 5, we introduce generalized hyperK\"ahler, $\G$ and Spin$(7)$ structures as special $\B(V)$-structures. The generalized exceptional structure (G$_2$ and Spin$(7)$ ) are discussed by Witt [W] from 
other point of view.
Our approach is adapted in these interesting cases. 
For instance, we will show that deformations of generalized Spin$(7)$structures are unobstructed. We will discuss the deformation problems of other special structures in a
forthcoming paper.
\head \S1. Clifford algebras and Spin groups \endhead
\subhead \S1-1\endsubhead 
Let $V$ be an $n$ dimensional real vector space and $V^*$ the dual space of $V$. 
We denote by $\eta(v)$ by the natural coupling between  $v\in V$ and $\eta\in V^*$. 
Then there is an indefinite bilinear form $\lan\,,\,\ran$ on the direct sum $\VV$
which is defined by 
$$
\lan E_1, E_2\ran =\frac12\eta_1(v_2)+\frac12\eta_2(v_1),
\tag1-1-1$$
where $E_i=v_i+\eta_i\in \VV$ for $i=1,2$.
(In particular the norm $\| E\|^2= \lan E, E\ran$.)
We consider $\VV$ as the $2n$ dimensional vector space and denote by 
$\otimes^k(\VV)$ the tensor product of $k$-copies of $\VV$. 
Let 
$$
\otimes(\VV):=\sum_{i=0}^\infty \otimes^k(\VV).
\tag1-1-2$$
be the the tensor algebra of $(\VV)$
(Note that $\otimes^0(\VV)=\R$),  and define
 $\I$ to be the two-sided ideal in $\otimes (\VV)$ generated by all elements of the form
 $E\otimes E-
\| E\| 1$ for $E\in \VV$.
Then the Clifford algebra $\CL(\VV)$ is defined to be the quotient algebra with the unit $1$ :
$$
\CL(\VV)=\otimes(\VV)/\I.
\tag1-1-3$$
The product of the Clifford algebra is called the Clifford product which is denoted by 
$\a\cdot\b$ for $\a,\b\in \CL(\VV)$ and for all $E, F\in \VV$,
$$
 E\cdot F+F\cdot E =\lan E, F\ran 1.
\tag1-1-4$$
Since the ideal $\I$ is generated by tensors of degree $2$, 
the Clifford algebra $\CL(\VV)$ is decomposed into the even part and the odd part : 
$$
\CL(\VV)=\CL^{\even}\oplus \CL^{\odd}.
\tag1-1-5$$
There are two involutions of $\CL(\VV)$ which play the smart roles.  
The first one is defined by the decomposition (1-1-5) : 
$$
\til\a :=
\cases  
&+\a,\quad  (\a\in \CL^{\even}),\\
&-\a, \quad (\a\in\CL^{\odd}),
\endcases
\tag1-1-6$$
for $\a\in \CL(\VV)$. 
If we reverse the order in a simple product  $\a=E_1\cdot E_2\cdots E_k\in \CL(\VV)$ ( $E_i\in \VV)$, 
we obtain the second involution $\sig$ of $\CL(\VV)$:
$$
\sig(\a)=E_n\cdots E_2\cdot E_1.
\tag1-1-7$$
Since there is the natural isomorphism between the skew-symmetric tensors $\w^*(\VV)$ and $\CL(\VV)$ as $\R$-module, 
there is the metric $\lan\,,\,\ran$ on $\CL(\VV)$ which is written as 
$$
\lan\a, \b\ran =\frac12\lan 1, \sig(\a)\b\ran,
\tag1-1-8$$
for $\a, \b\in \CL(\VV)$. 
In particular we denote by $\|\a\|$ the Clifford norm of $a$ : 
$$
\|\a\|^2:=\lan \a,\a\ran =\frac12\lan 1,\,\sig(\a)\a\ran.
\tag1-1-9$$
Let  $\w^pV^*$ be the space of skew-symmetric tensor of degree $p$ and $S$ the direct sum of the spaces of skew-symmetric tensors : 
$$
S:=\oplus_{p=0}^\infty \w^pV^*.
\tag1-1-10$$
Then $E=v+\eta\in \VV$ acts on $S$ by the interior and the exterior product : 
$$
E\cdot \phi =i_v\phi+\eta\w\phi
\tag1-1-11$$
Since we have the identity :
$$
i_v\eta\w\phi+\eta\w i_v\phi =\|E\|^2\phi,
\tag1-1-12$$
we have the action of $\CL(\VV)$ on $S$,(which is called the spin representation).
Let $\CL(\VV)^\times$ be the group which consists of invertible elements of $\CL(\VV)$.
For each $g\in \CL(\VV)^\times$ we define a linear map $\wtil\Ad_g \:\CL(\VV)\to \CL(\VV)$ by 
$$
 \wtil\Ad_g(\a) := \til g^{-1}\a g, \quad( \a\in \CL(\VV)),
\tag1-1-13$$
where $\til g$ is the first involution of $g$. 
Note that the image $\wtil\Ad_g(\VV)$ is not a subspace of $\VV$ for a general $g\in \CL(\VV)^\times$.
The conformal pin group $\CPin(=\CPin(\VV))$ is a subgroup of $\CL(\VV)^\times$ which defined by 
$$
\CPin:=\{\, g\in \CL(\VV)^\times\,|\,\wtil\Ad_g(\VV)\subset \VV\, \}.
\tag1-1-14$$ 
Since $\wtil\Ad_g$ is an orthogonal endmorphism of $\VV$,
we have the short exact sequence :
 $$\minCDarrowwidth{0.3cm} 
 \CD 
 1@>>>\R^\times @>>> \CPin@>\wtil\Ad>>\text{\rm O}(\VV)@>>>1.
 \endCD
 \tag1-1-15$$
 Since each element of the conformal pin group $\CPin$ is written as a simple product, it follows that
 the Clifford norm of $g\in \CPin$ is given by $\|g\|^2=\sig(g)\cdot g$.
 We define the Pin group $\Pin(=\Pin(\VV))$ by 
 $$
 \Pin=\{\, g\in \CPin\, |\, \|g\|=\pm1\, \},
\tag1-1-16$$
 and the Spin group $\Spin(=\Spin(\VV))$ is defined by 
 $$
 \Spin:=\Pin\cap\CL^{\even}.
 \tag1-1-17$$
 Then we also have the short exact sequence :
 $$
 \minCDarrowwidth{0.3 cm} 
 \CD
 1@>>>\Z_2@>>>\Spin@>\Ad>>\SO(\VV)@>>>1.
 \endCD
 \tag1-1-18$$
 We denote by $\Spin_0(=\Spin_0(\VV))$ the identity component of $\Spin$. 
 Then $\Spin_0$ is given by  
 $$
 \Spin_0=\{\, g\in \Spin\, |\, \|g\|=1\, \}.
 \tag1-1-19$$
\subhead \S1-2 \endsubhead
The Lie algebra so$(\VV)$ of the Lie group $\SO(\VV)$ is decomposed into three parts : 
$$
so(\VV)=\text{End} (V)\oplus \w^2V\oplus \w^2 V^*. 
\tag1-2-1$$
Each $a\in $so$(\VV)$ is written as a form of matrix : 
$$
\pmatrix 
A&\b\\
b&-A^*\endpmatrix,
$$
where $A\in $End$(V)$, $b\in \w^2V^*$, $\b\in \w^2V$ and $A^*\in$ End$(V)$ is defined by $A^*(\eta)(v)=\eta(Av)$ for  $v\in V$ and $\eta\in V^*$. 
(Note that $b\: V\to V^*$ and $\b\: V^*\to V$.) 
Then the Lie group $\GL(V)$ is embedded into $\SO(\VV)$ : 
$$
\pmatrix 
g&0\\
0&(g^*)^{-1}
\endpmatrix,\qquad  g\in \GL(V). 
\tag1-2-2$$
Further  
 for $b\in \w^2V^*$ and $\b\in \w^2V$ we  define $e^b$ and $e^\b$ by : 
 $$\align 
 &e^b=1+b+\frac1{2!}b^2+\cdots, \\
 &e^\b=1+\b+\frac1{2!}\b^2+\cdots,
 \endalign
 $$
 then $e^b$ and $e^\b$ are elements of Spin$_0$ respectively. 
 Let $\GL_0(V)$ be the identity component of $\GL(V)$. We denote by $q$ the embedding 
 (1-2-2) of $\GL_0(V)$ into the identity component of SO$_0(\VV)$ : 
 $$
 q\: \GL_0(V)\to \SO_0(\VV).
 \tag1-2-3$$
 Let $\Ad$ be the covering map $\Ad\: \Spin_0\to \SO_0(\VV)$ as before. 
 Then there is a map $p \: \GL_0(V)\to \Spin_0$ such that Ad$\circ p=q$, so that is, 
 $p$ is the lift of the map $q$ : 
 $$
 \minCDarrowwidth{0.3 cm} 
 \CD
 \GL_0(V)@>p>> \Spin_0 \\
@|  @VV\Ad V\\
 \GL_0(V)@>q>>\SO_0.
  \endCD
 \tag1-2-4$$
 The representation $S$ of the Clifford algebra $\CL(\VV)$ restricts to the representation $\rho_{\text{\rm spin}}$ of $\Spin_0$.
 $$
 \rho_{\text{\rm spin}}\: \Spin_0\to \GL(S).
 \tag1-2-5$$
 We also denote by $\rho_{\GL}^*$ the linear representation of $\GL(V)$ on $\w^*V^*$. 
 The composition $\rho_{\text{\rm spin}}\circ p$ gives rise to a representation of $\GL_0(V)$.
 \proclaim{Lemma 1-2-1} 
 The representation $\rho_{\text{\rm spin}}\circ p$ is given by 
 $$
 \rho_{\text{\rm spin}}=(\det V^*)^{\frac12}\otimes(\rho^*_{\GL})^{-1},
 $$
 where $(\det V^*)^{\frac12}$ is the half of the determinant representation.
 \endproclaim
  \head \S2. Clifford-Lie operators\endhead
We use the same notation as in section 1. Let $X$ be a real manifold of dim $n$. Then we consider the direct sum $T\oplus T^*$ of the tangent bundle $T=TX$ and the cotangent bundle $T^*=T^*X$. Let $\CL(X)=\CL(\TT)$ be the Clifford bundle on $X$ :
$$
\CL(X):=\bigcup_{x\in X}\CL(T_xX\oplus T^*_xX) \to X.
$$
We also define  the conformal pin group-bundle Cpin$(X)$=Cpin$(\TT)$ by : 
$$
\CPin(X):=\bigcup_{x\in X}\CPin(T_xX\oplus T^*_xX) \to X.
$$
Let $\pi$ be the natural 
projection, 
$$
\pi\: \otimes(\TT)\to \CL(X)=\otimes(\TT)/\I.
$$We denote by $\CL^{2i}$ the image 
$$
\CL^{2i}:=\pi\(\oplus_{l=0}^i\otimes^{2l}(\TT)\).
$$
Then we have a filtration of $\CL^{\even}$ :
$$
\CL^0\subset \CL^2\subset\CL^4\subset\cdots.
$$
We also have a filtration of $\CL^{\odd}$ which defined by
$$
\CL^1\subset\CL^3\subset\CL^5\subset\cdots,
$$
where
$$
\CL^{2i+1}:=\pi \(\oplus_{l=0}^i\otimes^{2l+1}(\TT)\).
$$
Let $S(X)$ be the bundle of differential forms $\w^*T^*X$
over a manifold $X$. 
By using the spin representation on each fibre as in section 1,
the bundle of the Clifford algebra $\CL(X)$ acts on $S(X)$.
Let $\L_{E}$ be the anti-commutator $\{ d, E\}=dE+Ed$ for a section $E$ of the bundle $\TT.$ 
(For simplicity, we denote it by $E\in \CL^1=\TT$.)
If we denote $E=v+\theta\in T\oplus T^*$
then $\L_{E}=\L_v+( d\theta)$, where $\L_v$ is the ordinary Lie derivative and $(d\t)$ acts on $S(X)$ by the wedge product.
Next we consider a bracket $[\L_E, F]=
\L_E F-F\L_E$ for $E,F\in T\oplus T^*$.
\proclaim{Lemma 2-1}
The bracket $[\L_E, F]$ is a section of
$T\oplus T^*$.
\endproclaim
\demo{proof}
When we write
$E=v+\theta$, $F=w+\eta\in T\oplus T^*$,
then we have
$$
\align
[\L_E, F]=&[\L_v+(d\theta), w+\eta]\\
=&[\L_v, w]+[\L_v,\eta]+[(d\theta), w]
+[(d\theta), \eta]\\
=&[v,w]+(\L_v\eta)+[(d\theta), w].
\endalign
$$
Since $[(d\theta), w]\in (T\oplus T^*)$,
we have the result.
\qed\enddemo
In this paper Clifford algebra valued Lie derivatives play an significant role.
\proclaim{Definition 2-2(Clifford-Lie operators)}
A Clifford-Lie operator of order $3$ on $X$ is a differential operator acting on $S(X)$ which is locally written as
$$
L=\sum_{i,j} a^{ij}E_i\L_{E_j}+ K,
$$
on every open set $U$ on $X$ for some $E_i\in \CL^1( TU\oplus T^*U)$,
$a_{ij}\in C^\infty(U)$ and 
$K\in \CL^3(TU\oplus T^*U). $
\endproclaim
Let $\{x_1,\cdots, x_n\}$ be a local coordinates of $X$. 
We denote by $v_i$ the vector field $\frac{\pa}{\pa x_i}$ and 
$\t^i=dx^i$.
Then the extrior derivative $d$ is locally written as
$$
d=\sum_{i=0}^n \t^i\w\L_{v_i}.
$$
Hence $d$  is the Clifford-Lie operator of order $3$.

Let $a$ be a section of $\CL^2(T\oplus T^*)$. Then we have
\proclaim{Lemma 2-3}
If $L$ is a Clifford-Lie operator of order $3$ then the commutator $[ L, a]$ is also
a Clifford-Lie operator of order $3$.
\endproclaim
\demo{proof}
Let $f$ be a function on $X$ and $E=v+\theta$ a section of $T\oplus T^*$.
Since we have
$$
\L_{E}fa= (\L_E f)a+ f\L_Ea,
$$
where $\L_E f=\L_vf\in C^\infty(X)$.
We have
the following equality on an open set $U$ on $X$:
$$\align
[L, fa]=&L(fa)-faL\\
=&\sum_{ij}a_{ij}E_i\L_{E_j}(fa)-faL+K\\
=&\sum_{ij}a_{ij}E_i(\L_{E_j}f)a+f[L,a].
\endalign
$$
Since $E_i(\L_{E_j}f)a\in\CL^3(T\oplus T^*)$, it is sufficient to show the lemma in the case $a=F_1F_2$ for $F_i\in T\oplus T^*(i=1,2)$.
The bracket $[\L_E, F_1F_2]$ is given by
$$
\align
[\L_E, F_1F_2]=&
\L_E F_1F_2- F_1F_2\L_E\\
=&[\L_E, F_1]F_2+F_1\L_EF_2-F_1F_2\L_E\\
=&[\L_E, F_1]F_2+ F_1[\L_E, F_2].\\
\endalign
$$
Hence it follows from lemma 2-1 that $[\L_E,F_1F_2]\in \CL^2$.
The bracket $[E_1\L_{E_2}, F_1F_2]$ is given by
$$\align
[E_1\L_{E_2}, F_1F_2]=&E_1\L_{E_2}F_1F_2
-F_1F_2E_1\L_{E_2}\\
=&E_1[\L_{E_2}, F_1F_2]+
E_1F_1F_2\L_{E_2}-F_1F_2E_1\L_{E_2},\\
=&[E_1,\,F_1F_2]\L_{E_2} +E_1[\L_{E_2}, F_1F_2].
\endalign
$$
Since $[E_1, F_1F_2]=2\lan E_1 F_1\ran F_2-2\lan E_1, F_2\ran F_1\in \CL^1=(T\oplus T^*)$, it follows that
the bracket $[E_1\L_{E_2}, F_1F_1]$ is
a Clifford-Lie operator of order $3$.
Then the result follows from the equation:
$$\align
[L, F_1F_2]=&[\sum_{ij}a_{ij}E_i\L_{E_j},
F_1F_2] \\
=&\sum_{ij}a_{ij}[E_i\L_{E_j}, F_1F_2].
\endalign
$$
\qed\enddemo
\proclaim{Lemma 2-4}
The commutator $[ d, a]$ is a Clifford-Lie operator of order $3$.
\endproclaim
\demo{proof}
Since $d$ is a Clifford-Lie operator of order $3$, the result follows from lemma 2-3.
We think a following direct proof is more readable.
$[d, fa]=dfa-fad=(df)a+f[d,a]$ for a function $f$. Hence it is sufficient to show the lemma in the case $a=E_1E_2$, where $E_i\in T \oplus T^*$ $(i=1,2)$.
Then the bracket $[d,a]$ is written as
$$
\align
[d,a]=&dE_1E_2-E_1E_2d\\
=&\L_{E_1}E_2-E_1dE_2-E_1E_2d\\
=&\L_{E_1}E_2-E_1\L_{E_2}\\
=&E_2\L_{E_1}-E_1\L_{E_2}+[\L_{E_1},E_2].
\endalign
$$
Hence the result follows from $[\L_{E_1},E_2]\in \CL^1\subset\CL^3$.
\qed\enddemo
\proclaim{Proposition 2-5}
Let $a_1,a_2\in $CL$^2(T\oplus T^*)$.
Then $[[d,a_1], a_2]$ is a Clifford-Lie operator of order $3$.
Further we denote by Ad$_{a}L$ the commutator $[L, a ]$.
Then the composition
$Ad_{a_1}(Ad_{a_2}(\cdots Ad_{a_n} d)\cdots )$ is a Clifford-Lie operator of order $3$ for $a_1,\cdots , a_n\in \CL^2$.
\endproclaim
\demo{proof}
It follows form Lemma 2-3 and 4.
\qed\enddemo
\proclaim{Remark 2-6}
In the case of $a_1, a_2\in $End$(TX)$, the bracket 
$[[d,a_1],a_2]$ is given in terms of the Nijenhuis tensor of $a_1$ and $a_2$. 
In the case $a_1, a_2\in \w^2T$, the bracket 
$[[d,a_1],a_2]$ is the Schouten bracket. 
In general the bracket $[[d,a_1],a_2]$ is not a tensor but a differential operator.
\endproclaim
Let $a$ be a section of $\CL^2$ and $L$ an operator acting on $S(X)$. We successively define an operator $(\Ad_a^l)L$ acting on $S(X)$ by 
$$
(\Ad_a)^lL=[(\Ad_a)^{l-1}L, a].
$$
We also define a formal power series $\(\exp(\Ad_a)\)L$ by 
$$\align
\(\exp(\Ad_a)\)L=&\sum_{l=0}^\infty \frac1{l!}(\Ad_a)^lL\\
=&d+[L,a]+\frac1{2!}[[L,a],a]+\cdots.
\endalign
$$
\proclaim{Lemma 2-7} The power series $\(\exp(\Ad_a)\)L$ is given by 
$$
\(\exp(\Ad_a)\)L=e^{-a} \circ L\circ e^{a}.
$$
\endproclaim
\demo{proof} 
It follows from definition of $(\Ad_a)^l L$ that 
$$
(\Ad_a)^lL =\sum_{m=0}^l\frac{(-1)^m l!}{m!(l-m)!}a^m\, L\, a^{l-m}.
$$
Then by a combinatorial calculation we have 
$$
L\,a^k =\sum_{l=0}^k\frac{k!}{l!(k-l)!}\,a^{k-l}\,(\Ad_a)^lL.
$$
Then we have 
$$\align
L\, e^a=&e^a(L+(\Ad_a)L+\frac1{2!}(\Ad_a)^2L+\frac1{3!}(\Ad_a)^Ld+\cdots)\\
=&e^a(\exp(\Ad_a) L).
\endalign
$$
Hence the result follows.
\qed\enddemo
\proclaim{Proposition 2-8} 
If $L$ is a Clifford-Lie operator of order $3$ and $a\in \CL^2$, then 
$e^{-a}\circ L\circ e^a=\(\exp(\Ad_a) L\)$ is also a Clifford-Lie operator of order $3$. 
In particular $\(\exp(\Ad_a)d\)$ is a Clifford-Lie operator of order $3$.
\endproclaim
 \demo{proof} 
The result follows from lemmas 2-5 and 2-7.
 \qed\enddemo
 \head \S3. Deformations of generalized geometric structures\endhead 
\subhead \S3-1\endsubhead
Let $V$ be an $n$ dimensional real vector field and 
$V^*$ the dual space of $V$. As in section 1 the space of the skew-symmetric tensors S:=$\w^*V^*$ is regarded as the spin representation of $\CL(=\CL(\VV))$, which restricted to give the representation of 
$\CPin(=\CPin(\VV))$.
We consider the direct sum of the spin representations of $\CPin(\VV)$ :
$$\oplus^lS:= \overset{l\text{\rm times} }\to{\overbrace{(\w^*V^*\oplus\cdots\oplus\w^*V^*)}}.
$$
Let $\Phi_V=(\phi_1,\cdots, \phi_l)$ be an element of the direct sum 
$\oplus^l S$. Then we have the orbit $\B(V)$ of CPin$(\VV)$ through $\Phi_V$ :
$$
\B(V):=\{\, g\cdot\Phi_V\, |\, g\in \CPin(\VV)\,\}.
$$
From now on we fix the orbit $\B(V)$.
In this paper we think that 
the orbit of $\CPin(\VV)$ gives rise to a generalized geometric structure on the vector space. 
Let GL$_0(V)$ be the connected component of GL$(V)$ with the identity
and $\A(V)$ the orbit of  GL$_0(V)$ through $\Phi_V$. 
As in the section 1, we have the map $p\:$ GL$_0(V)\to \Spin_0\subset\CPin$.
It follows from lemma 1-2-1 that for $\Phi\in \oplus^lS$ we have,
$$
(p(g))\cdot\Phi=(\det g)^{\frac12}\rho_{\GL^*(V)}(g^{-1})\Phi,\quad \text{\rm for } g\in \GL_0(V). 
$$ 
Since $(\det g)^{-\frac12}p(g)\in \CPin$, we have 
$$
((\det g)^{-\frac12}p(g))\cdot\Phi=\rho_{\GL^*(V)}(g^{-1})\Phi.
$$
It implies that the GL$_0(V)-$orbit $\A(V)$ is embedded into the CPin$(\VV)-$orbit $\B(V)$ : 
$$
\A(V)\hookrightarrow \B(V).
\tag3-1-1$$
The inclusion (3-1-1) shows that the group $\CPin$ is suitable for our construction.
Let $X$ be an oriented and compact real manifold of dim $n$. 
As in section 2 we have the Clifford bundle $\CL(X)$ and the conformal pin bundle 
$\CPin(X)$ over $X$.
When we take an identification between $V$ and $T_xX$ for each $x\in X$, 
we have the orbit $\B(T_xX)$ of CPin$(T_xX\oplus T^*_xX)$. It follows from (3-1-1) that the orbit $\B(T_xX)$ is independent of a choice of an identification and thus $\B(T_xX)$ is canonically defined as the submanifold of the direct sum of forms $\oplus^l \w^*T^*_xX$. Hence we have the fibre bundle $\B(X)\to X$: 
$$
\B(X):=\bigcup_{x\in X}\B(T_xX)\to X.
$$
Let $H$ be the isotropy group of the action of CPin$(\VV)$ at $\Phi_V$:
$$
H:=\{\, g\in \CPin(\VV)\, |\, g\cdot\Phi=\Phi\, \}.
$$
Then the fibre bundle $\B(X)$ is the fibre bundle with fibre 
$\CPin(\VV)/ H$ and $\B(X)$ is embedded into the direct sum of differential forms $\oplus^l \w^*T^*X$. We denote by $\E_B(X)$the set of C$^\infty$-sections of the fibre bundle $\B(X)$: 
$$
\E_B(X):=C^\infty(X, \B(X)).
$$
Each section $\Phi\in \E_B(X)$ consists of differential forms on which the exterior derivative $d$ acts. 
Let $\wtil\M_B(X)$ be the set of $d$-closed section of $\E_B(X)$: 
$$
\wtil\M_B(X):=\{\, \Phi\in \E_B(X)\, |\, d\Phi=0\, \}.
$$
\proclaim{Definition 3-1} 
A generalized geometric structure on $X$ associated with the orbit $\B(V)$ 
is a $d$-closed section $\Phi\in \wtil\M_B(X)$. 
For simplicity, we call a $d$-closed section $\Phi$ a $\B(V)$-structure on $X$.  
\endproclaim
The diffeomorphism $\Diff(X)$ naturally acts on $\wtil\M_B(X)$ by the pull back. We denote by $\Diff_0(X)$ the identity component of $\Diff(X)$.
Since the exponential $e^{d\gam}$ is a section of the bundle Spin$_0(X)$ for a $1$-form $\gam$, we have the action of $e^{d\gam}$ on $\B(V)$-structures $\wtil\M_{\B}(X)$,
$$
\Phi \mapsto e^{d\gam}\w\Phi,\qquad ( \gam \in T^*X).
$$
Let $\wtil\Diff_0(X)$ be the group generated by the composition of the action of $\Diff_0(X)$ and $d$-exact $2$-forms: 
$$
\wtil{\Diff_0}(X):= \{\, e^{d\gam}\w f^*\, |\, \gam\in T^*,\,f\in \Diff_0(X)\, \}.
$$
Here the group $\wtil\Diff_0(X)$ is regarded as a subgroup of the automorphisms of the bundle Spin$_0(X)$ : 
$$
\CD 
\Spin_0(X)@>>>\Spin_0(X)\\
@VVV @VVV\\
X@>>>X.
\endCD
$$
Hence the group $\wtil\Diff_0(X)$ is an extension of $\Diff_0(X)$ by $d$-exact $2$-forms $d\,(\w^1 T^*$): 
$$\minCDarrowwidth{0.3cm} 
\CD
0@>>>d\,(\w^1T^*)@>>>\wtil\Diff_0(X)@>>>\Diff_0(X)@>>>0.
\endCD
$$
\proclaim{Definition 3-2} 
A moduli space $\M_B(X)$ of $\B(V)$-structures on $X$ is the quotient space of $\wtil\M_B(X)$ divided by the action of $\wtil\Diff_0(X)$: 
$$
\M_B(X):= \wtil\M_B(X)/\Diff_0(X).
$$
\endproclaim
\subhead \S 3-2 \endsubhead
Let $\B(V)$ be the fixed orbit of CPin$(\VV)$ as in section 3-1 and $\Phi$ a $\B(V)$-structure on a manifold $X$. 
In order to consider deformations of $\Phi$, we introduce a deformation complex of the $\B(V)$-structure $\Phi$. 
As in section 2 there is the natural filtration of the even Clifford bundle CL$^{\even}$ and 
the one of the odd Clifford bundle $\CL^{\odd}$)  :
$$\align
&\CL^0\subset \CL^2\subset\CL^4\subset\cdots, \\
&\CL^1\subset\CL^3\subset\CL^5\subset\cdots.
\endalign 
$$
Then by using the action of $\CL^k$ on $\Phi$, we obtain a vector bundle $\bold E^{k-1}(X)$ over $X$: 
$$
\bold E^{k-1}(X):=\CL^k\cdot\Phi,
$$
and the corresponding filtrations of vector bundles: 
$$\align
&\bold E^{-1}(X)\subset \bold E^1(X)\subset\bold E^3(X)\subset\cdots,\\
&\bold E^0(X)\subset\bold E^2(X)\subset\bold E^4(X)\subset \cdots.
\endalign
$$
(Note that we shift the degree of the filtration of vector bundles.)
The vector bundle $\bold E^{-1}(X)$ is the line bundle generated by $\Phi$. The vector bundle $\bold E^0(X)$ is generated by $E\cdot \phi$ 
for all $E\in \TT$ over C$^\infty(X)$ and $\bold E^1(X)$ is generated by 
$E_1\cdot E_2\cdot \Phi$ for all $E_1, E_2\in \TT$. 
Each $\bold E^k(X)$ is embedded into the direct sum of differential forms on which the exterior derivative $d$ acts. 
\proclaim{Proposition 3-3} There is a differential complex $\#_\B(=\#_{\B,\Phi})$ for each $\Phi\in\wtil\M_{\B}(X)$, 
$$\minCDarrowwidth{0.2cm} 
\CD 
0@>>>\bold E^{-1}(X)@>d_{-1}>>\bold E^0(X)@>d_0>>\bold E^1(X)@>d_1>>
\bold E^2(X)@>d_2>>\cdots,
\endCD
$$
where $d_k$ is given by the restriction $ d|_{\bold E^k(X)}$. The cohomology groups of the complex $\#_\B$ is denoted by H$^k(\#_\B)$,
$$
H^k(\#_\B):=\frac{\ker d_k \:\Gam( \bold E^k(X))\to\Gam(\bold E^{k+1}(X))}
{\text{\rm im}\, d_{k-1}\:\Gam(\bold E^{k-1}(X))\to \Gam(\bold E^k(X))}.
$$ 
Then the first cohomology group $H^1(\#_\B)$ is regarded as the infinitesimal tangent space of the deformations of the $\B(V)$-structure $\Phi$.
\endproclaim
\demo{proof} 
A section of $\bold E^{-1}(X)$ is written as $f\Phi$ for a function $f$. 
Hence $d( f\Phi)= df\w \Phi$ and we see that the image $d(\bold E^{-1}(X))$
is included into $\bold E^0(X)$. We denote by $\L_{F}$ the anti-commutator $d F+Fd$ acting on forms where $F\in \TT$. 
When we write $F=v+\eta$ for $v\in T$ and $\eta\in T^*$, then $\L_F$ is given by 
$$
\L_F=\L_v+(d\eta)\w,
$$
where $\L_v$ denotes the Lie derivative. Then we have 
$$\align
\L_F(f\Phi)=&\L_v(f\Phi)+(d\eta)\w(f\Phi)\\
=&(\L_vf)\Phi+f\L_v\Phi+f(d\eta)\w\Phi,
\endalign
$$
where $\L_vf\in C^\infty(X)$.
Since $\GL_0(TX)$ is the subbundle of $\CPin(X)$, diffeomorphisms of $X$ acts on $\E_\B(X)$. Hence we have 
$$\L_v\Phi\in T_\Phi\E_\B(X).
$$
The conformal Spin$_0$ bundle $\CPin_0(X)$ is given by 
$$
\CPin_0(X)=\{\, e^a\, |\, a\in \CL^2\,\}.
$$
Since the tangent space $T_\Phi\E_\B(X)$ is generated by the action of $\CPin_0(\TT)$, 
we have 
$$
T_\Phi\E_\B(X)\cong\CL^2\cdot \Phi=\bE^1(X).
$$
Hence we have 
$$
\L_v\Phi\in \bE^1(X).
$$
Then it follows that $\L_F(\bE^{-1}(X))\subset \bE^1(X)$.
We also have 
$$
d(F\cdot\Phi)=\L_F\Phi-Fd\Phi=\L_F\Phi.
$$
Hence we have $d(\bE^0(X))\subset \bE^1(X)$.
For $F_1, F_2\in \TT$ we have 
$$
\L_{F_1}(F_2\cdot\Phi)=[\L_{F_1}, F_2]\Phi+F_2\cdot\L_{F_1}\Phi.
$$
It follows from lemma 2-1 that $[\L_{F_1}, F_2]\in \TT$. 
Hence from $\L_{F_1}\Phi\in \bE^1(X)$ we have 
$\L_F(\bE^0(X))\subset\bE^2(X)$. We will show that $d\bE^k(X)\subset \bE^{k+1}(X)$ by induction on $k$.
We assume that 
$d\bE^{k-2}(X)\subset \bE^{k-1}(X)$ and $\L_F(\bE^{k-2}(X))\subset \bE^k(X)$ for some 
$k\geq 1$ and for all $F\in \TT$. 
Then for $F_1, F_2\in \TT$ and $s\in \bE^{k-2}(X)$ we have 
$$\align
d(F_1\cdot F_2\cdot s)=&\L_{F_1}(F_2\cdot s)-F_1\cdot dF_2\cdot s\\
=&[\L_{F_1}, F_2]\cdot s+F_2\cdot\L_{F_1}s\\
&-F_1\cdot\L_{F_2} s+F_1\cdot F_2 \cdot ds.
\endalign
$$
It follows from our assumption ($ds\in \bE^{k-1}(X)$ and $\L_F s\in \bE^k(X)$) that 
$d(F_1\cdot F_2\cdot s)\in \bE^{k+1}(X)$ since $[\L_{F_1}, F_2]\cdot s\in \bE^{k-1}(X)
\subset\bE^{k+1}(X)$.  
Hence $d(\bE^k(X))\subset \bE^{k+1}(X)$.
For $F_3\in \TT$ we also have 
$$\align
\L_{F_3}(F_1\cdot F_2\cdot s)=&[\L_{F_3}, F_2]\cdot F_1\cdot s+F_2\cdot\L_{F_3}(F_1\cdot s),\\
=& [\L_{F_3}, F_2]\cdot F_1\cdot s+F_2\cdot[\L_{F_3},F_1] \cdot s\\
&+F_2\cdot F_1\cdot \L_{F_3} s.
\endalign
$$
Hence it follows from our assumption $\L_F s\in \bE^k(X)$ that 
$\L_{F_3}(F_1\cdot F_2\cdot s)\in \bE^{k+2}(X)$. 
Hence $\L_F(\bE^k(X))\subset\bE^{k+2}(X)$. 
We already show that our assumption in cases of $k=1,2$.
Therefore we have $d\bE^k(X)\subset\bE^{k+1}(X)$ for all $k$ by induction.
The tangent space of the orbit of $\wtil\Diff_0(X)$ is given by the Lie derivative $\L_v\Phi$ and $d\gam\w\Phi$ for $v\in T$ and $\gam \in T^*$.
Hence it follows that the image $d(\Gam(E^0(X))$ is the tangent space of $\wtil\Diff_0(X)$. 
As we see, the tangent space of $\E_\B(X)$ is global sections of $\bold E^1(X)$. 
Hence the infinitesimal tangent space of deformations of $\Phi$ is given by 
the first cohomology group $H^1(\#_\B)$.
\qed\enddemo
The direct sum $\oplus^l S(=\oplus^l\w^*T^*)$ is invariant under the action 
of the exterior derivative $d$ which is the direct sum of the full de Rham complex. For simplicity we call $\oplus^l S$ the de Rham complex.
Then the complex $\#_\B$ is the subcomplex of the de Rham complex:
$$
\minCDarrowwidth{0.3cm}
\CD
0@>>>\bE^{-1}(X)@>d_{-1}>>\bE^0(X)@>d_0>>\bE^1(X)@>d_1>>\bE^2(X)@>>>\cdots,\\
@. @VVV @VVV @VVV @VVV \\
\cdots @>>>\oplus^l \w^*T^*@>d>>\oplus^l \w^*T^*@>d>>\oplus^l 
\w^*T^*@>d>> \oplus^l \w^*T^*@>>>\cdots
\endCD
$$
We denote by $H^*_{dR}(X)(=\oplus^l \oplus_{p=0}^{\dim X} H^p(X,\R) )$ the cohomology group of the de Rham complex. Then we have the map $p^k_\B$ : 
$$p^k_\B\:  H^k(\#_\B)\to H^*_{dR}(X).
$$
Since the action of $\wtil\Diff_0(X)$ on $\wtil\M_\B(X)$ preserves a de Rham cohomology class $[\Phi]$ of $\B(V)$-structure $\Phi$, we have the map $P_\B$: 
$$
P_\B\: \M_\B(X) \to H^*_{dR}(X).
$$
The map $P_\B$ is called the period map.
\proclaim{Definition 3-4} 
An orbit $\B(V)$ is completely elliptic if the differential complex $\#_\B$ is a elliptic complex. 
In particular, an orbit $\B(V)$ is elliptic if the complex is elliptic at degrees $k=1,2$. 
\endproclaim
\proclaim{Definition 3-5} Let $\B(V)$ be an orbit of $\CPin\,(\VV)$ as in before. 
We say a $\B(V)$-structure $\Phi$ on $X$ is topological if the map $p^k_\B\:H^k(\#_\B)\to H^*_{dR}(X)$ is injective for $k=1,2$. 
An orbit $\B(V)$ is topological if each $\B(V)$-structure $\Phi$  is topological on every compact and oriented $n$-manifold.
\endproclaim
The complex $\#_\B$ is elliptic if the corresponding symbol complex is exact. Hence the elliptic condition only depends on a choice of an orbit $\B(V)$.
However the topological condition is depending on a choice of a $\B(V)$-structure $\Phi$ on $X$.
\proclaim{Definition 3-6} 
A $\B(V)$ structure $\Phi$ on $X$ is unobstructed if for each representative element $\a$ of the infinitesimal tangent space $\H^1(\#_\B)$, there exists one parameter family of deformations $\Phi_t\in \wtil\M_\B(X)$ with $\Phi_0=\Phi$ such that 
$$
\frac{d}{dt}\Phi_t|_{t=0}=a. 
$$
\endproclaim
If $\Phi$ is unobstructed, each infinitesimal tangent generates an actual deformations and 
the space of deformations of $\Phi$ is locally given by an open set of $H^1(\#_\B)$.
 From the view point as in [Go] , we have the following criterion for unobstructed deformations of $\B(V)$-structures and the Torelli-type theorem :
 \proclaim{Theorem 3-7} 
 Let $\B(V)$ be an elliptic orbit and
 $\Phi$ a $\B(V)$-structure on a compact and oriented $n$-manifold $X$.
 If $\Phi$ is topological, then deformations of $\Phi$ are unobstructed and  the deformations of $\Phi$ is locally embedded into the de Rham cohomology group $H^*_{dR}(X)$.
 In particular, if an orbit $\B(V)$ is elliptic and topological, the period map $P_\B$ of the moduli space 
 $\M_\B(X)$ of $\B(V)$ structures on $X$ is locally injective.
 \endproclaim
 Our proof of theorem 3-7 is shown by the following theorems 3-8 and 3-9. 
\proclaim{Theorem 3-8} 
Let $\B(V)$ be  an elliptic orbit of $\CPin(\VV)$ and $\Phi$ a $\B(V)$-structures on a compact, oriented manifold $X$. 
If the map $p^1_\B$ for $\Phi$ is injective, then there exists a neighborhood $U$ of $\Phi$ in the moduli space $\M_\B(X)$ such that 
the restriction of the period map $P_\B|_U\: U\to H^*_{dR}(X)$ is injective. 
\endproclaim
(Note that Proposition is regarded as a generalization of the Moser's stability theorem for symplectic structures an volume forms.)
\proclaim{Theorem 3-9}
Let $\B(V)$ be  an elliptic orbit of $\CPin(\VV)$ and $\Phi$ a $\B(V)$-structure on a compact, oriented manifold $X$. 
If $p^2_\B$ for $\Phi$ is injective then deformations of $\Phi$ are unobstructed.
\endproclaim
Our proof of theorem 3-9 is shown in next section 3-3. 
In order to obtain theorem 3-8, we will show the following lemma 
\proclaim{Lemma 3-10} 
Let $\{\Phi_n\}_{n=1}^\infty$ be a sequence of $\B(V)$-structures which converges to a $\B(V)$ structure $\Phi$, so that is, 
$$
\lim_{n\to\infty}\Phi_n =\Phi\in \wtil\M_\B(X).
$$
We denote by $\bE^k_n(X)$ the vector bundle which defined by $\bE^k_n(X)=\CL^{k+1}\cdot\Phi_n$ and $\#_{\B, n}$ by the deformation complex $\{\bE^*_n\}$ with 
cohomology groups $H^k(\#_{\B,n})$. If the map $p^1_{\B, n}\: H^1(\#_{\B,n}) \to H^*_{dR}(X)$ is not injective for all $n$, then 
the map $p^1_\B$ with respect to $\Phi$ is also not injective.
\endproclaim
Lemma 3-10 shows that the injectivity of the map $p^1_\B$ is an open condition, so that is, 
if $p^1_\B$ is injective for $\Phi\in \wtil\M_\B(X)$, then there exists an neighborhood $\til U$ such that $p^1_\B$  is also injective for all $\Psi\in\til U$.
\demo{proof of lemma 3-10} We take a Riemannian metric on the manifold $X$.
Then we have  the Laplacian  $\triangle_{\B, n}=d_1^*d_1+d_0d_0^*$  defined by the complex $\{\bE^*_n\}$ acting on sections of $\bE^1_n(X)$. We denote by $\Bbb H^1(\#_{\B,n})$ the kernel of the Laplacian $\triangle_{\B,n}$.
Since the complex $\#_{\B,n}$ is elliptic, the cohomology group $H^1(\#_{\B,n})$ is isomorphic to $\Bbb H^1(\#_{\B,n})$. We also have the ordinary Laplacian $\triangle$ which acts on $\oplus^lS$ and we denote by $\Pi$ the $L^2$-projection to the $\trian$-Harmonic forms. If $p^1_{\B,n}$ is not injective, we have $a_n\in \CL^2$ such that 
$a_n\cdot\Phi_n$ is a non-zero element of  $\Bbb H^1(\#_{\B,n})$ with $\Pi(a_n\cdot\Phi_n)=0$. 
For each $\Phi_n$ we can take  a section $g_n$ of the fibre bundle $\CSpin_0(X)$ with 
$g_n\cdot \Phi_n =\Phi$ and $g_n\to 1$ as $n\to \infty$.
By the left multiplication $L_{g_n}$ of $g_n$, we identify $\bE^1_n(X)$ with $\bE^1(X)=\CL^2(X)\cdot\Phi$, 
$$\align
L_{g_n}\: &\bE^1_n(X) \to \bE^1(X),\\
&a_n\cdot\Phi_n \mapsto g_n\cdot a_n\cdot\Phi_n=(\Ad_{g_n}a_n )\cdot\Phi.
\endalign
$$
Then the elliptic operator $\wtil\trian_{\B,n}$ on $\bE^1(X)$ is induced by 
$$
\wtil\trian_{\B,n}=L_{g_n}\trian_{\B,n}L_{g_n}^{-1}.
$$
We put $b_n=\Ad_{g_n}a_n$. 
Then we have 
$$
\wtil\trian_{\B,n}b_n\cdot\Phi=L_g\trian_{\B,n}(a_n\cdot\Phi_n)=0.
$$
We take $a_n$ such that the Sobolev norm of $b_n\cdot \Phi$ is normalized, 
$$
\|b_n\cdot\Phi\|_{L_4^2}=1.
$$
Then from Rellich lemma there exists a subsequence $\{ b_m\cdot\Phi\}_m$ which converges to $b\cdot\Phi\in \bE^1(X)$ with respect to the norm $L_2^2$.
Since $\wtil\trian_{\B,m}b_m\cdot\Phi=0$, we have an estimate, 
$$
\|b_m\cdot\Phi\|_{L_4^2}\leq C_1\|b_m\cdot\Phi\|_{L^1}\leq C_2\|b_m\cdot\Phi\|_{L_2^2},
$$
where $C_i\neq 0$ does not depend on $m$ for $i=1,2$.
Hence we have the bound,
$$
0 \neq C_3\leq \|b\cdot\Phi\|_{L_2^2}.
$$
The family of elliptic operator $\{\wtil\trian_{\B,m}\}_m$ also converges to the operator
$\trian_\B$ as $m\to \infty$. Hence we have 
$$
\trian_\B(b\cdot\Phi)=0.
$$
Since $g_m\to 1$($m\to \infty)$, the sequence $\{a_m\cdot\Phi_m\}=\{g_m^{-1}\cdot
b_m\cdot\Phi\}_m$ converges to $b\cdot\Phi$ ($n\to \infty$).  Hence  it follows from
$\Pi(a_m\cdot\Phi_m)=0$,
$$
\Pi (b\cdot\Phi)=0.
$$
Hence $b\cdot\Phi\neq 0$ is an element of $\ker p^1_\B$ and we have the result.
\qed \enddemo
\demo{proof of theorem 3-8} 
Let $\til U$ be a neighborhood of $\Phi$ such that $p^1_{\B}$  is injective for every $\Psi\in \til U$. 
Let $\{\Phi_t\}$, $0\leq t\leq 1$ be a smooth family of $\B(V)$-structures in the neighborhood $\til U$. 
We assume that the $d$-closed form $\Phi_t$ belongs to the same de Rham cohomology class as $\Phi_0$ for all $t$, so that is, there exists $A_t$ such that 
$$
\Phi_t-\Phi_0=dA_t.
\tag3-2-1$$
Since the group $\wtil\Diff_0(X)$ is generated by the action of $\Diff_0(X)$ and the action of $d$-exact $b$-fields,  theorem is reduced to the followings : 
\enddemo
\proclaim{Proposition 3-11} If the map $p^1_\B$ is injective for all $\Phi_t$, then there exist a smooth family of differmorphisms $\{f_t\}$ and a smooth family of $d$-exact $2$-forms $\{d\gam_t\}$ such that 
$$
e^{d\gam_t}\w f_t^*\Phi_t=\Phi_0,\quad \text{\rm for all } t\in [0,1].
\tag3-2-2$$
\endproclaim
\demo{proof of proposition 3-11}
By differentiating the equation (3-2-2), we have 
$$
\frac{d}{dt}\( e^{d\gam_t}\w f_t^*\Phi_t\) =0, \quad \forall \,t\in [0,1],
\tag3-2-3$$
which is equivalent to 
$$
e^{d\gam_t}\w d\dot\gam_t\w f_t^*\Phi_t+e^{d\gam_t}\w \dot{f}_t^*\Phi_t+e^{d\gam_t}\w f_t^* \dot{\Phi_t}=0.
\tag3-2-4$$
By the left action of $(f_t^{-1})^* (e^{-d\gam_t})$, we have
$$
(f_t^{-1})^*( d\dot\gam_t \w f_t^*\Phi_t)+(f_t^{-1})^*f_t^*\Phi_t+\dot{\Phi_t}=0.
\tag3-2-5$$
We set $ (f_t^{-1})^*\dot\gam_t=\til\gam_t$.
Since $(f_t^{-1})^*\dot{f_t}^*\Phi_t$ is given as the Lie derivative $\L_{v_t}\Phi_t$ for a vector field $v_t$, it follows from (3-2-1) that
$$
(d\til\gam_t)\w\Phi_t+\L_{v_t}\Phi_t+d\dot{A_t}=0.
\tag3-2-5$$
Since $\Phi_t$ is $d$-closed, we have 
$$
\L_{v_t}\Phi= di_{v_t}\Phi_t.
\tag3-2-6$$
We substitute (3-2-6) in (3-2-5) and we have 
$$\align
&(
(d\til\gam_t)\w\Phi_t+d i_{v_t}\Phi_t )=d\,\,(\til\gam_t+v_t)\cdot\Phi_t =-d\dot{A_t},
\tag3-2-7\endalign
$$
where $(v_t+\til\gam)\in \TT$ acts on $\Phi$ by the Clifford multiplication. 
We denote by $\bE^k_t(X)$ the vector bundle $\CL^{k+1}\cdot\Phi_t$ and 
$\#_{\B,t}$ the complex $\{\bE^*_t(X)\}$
We denote by $\bE^k_t(X)$ the vector bundle $\CL^{k+1}\cdot \Phi_t$ and 
by $\#_{\B, t}$ the complex $(\bE^*_t(X), d)$.
Then $(\til\gam_t+v_t)\cdot\Phi$ is a section of $\bE^0_t(X)$ and $-\dot{\Phi}_t=-d\dot A_t$ is 
a section of $\bE^1_t(X)$. Hence $-d\dot{A_t}$ yields the class $-[d\dot{A_t}]\in H^1(\#_{\B,t})$ of the deformation complex $\#_{\B,t}$ :
$$
\minCDarrowwidth{0.3cm} 
\CD
\bE_t^0@>d_0>>\bE_t^1@>d_1>>\cdots.
\endCD
$$
Then we see that the class $[-d\dot{A_t}]\in H^1(\#_{\B,t})$ vanishes
since the class $-[d\dot{A_t}]$ is represented by the $d$-exact form and 
the map $p^1_{\B,t}$ is injective. 
If we take a metric on the manifold $X$, we have the adjoint operator $d_k^*$ and the 
Green operator $G_t$ of the complex $\#_{\B,t}$. 
We define a section $B_t$ of $\bE^0_t(X)$  by 
$$
B_t=-d_0^*G_td\dot{A_t}.
\tag3-2-8$$
Then from the Hodge theory of the elliptic complex, we have 
$$
dB_t=-d\dot{A_t}.
\tag3-2-9$$
Since $B_t$ is written as $E_t\cdot\Phi$ for $E_t\in \TT$, we set $v_t+\til\gam_t=E_t$ such that a smooth family $\{v_t+\til\gam_t\}$ satisfies the equation (3-2-7). 
By solving the equation $(f_t^{-1})^*\dot{f_t}^*\Phi_t=\L_{v_t}\Phi$, we have the smooth family $\{f_t\}$ with $f_0=$id. 
Hence we have $\{f_t\}$ and $\{d\gam_t\}$ which satisfy the equation (3-2-2)
\qed\enddemo
\subhead 
\S 3-3  Construction of deformations
\endsubhead
This subsection is devoted to proof of theorem 3-9.
\demo{proof of theorem 3-9} Let $X$ be an $n$-dimensional, compact and oriented manifold with a $\B(V)$-structure $\Phi$. We take a Riemannian metric on $X$. 
(Note that this metric is independent to the structure $\Phi$.) 
The conformal pin bundle $\CPin(X)=\CPin(\TT)$ acts on the fibre bundle $\B(X)$ transitively. 
Hence every global section $\E_\B(X)$ is written as $g\cdot\Phi$ for a section $g$ of $\CPin(\TT)$. 
The identity component $\CSpin_0(\TT)$ of $\CPin(\TT)$ is given by 
$$
\CSpin_0(\TT)=\l\{\, e^a\, |\, a\in \CL^2(\TT)\,\r\}.
\tag3-3-1$$
Hence every deformation of $\Phi$ in $\E_\B(X)$ is given by 
$e^a\cdot \Phi$ for a section $a$ of $\CL^2(\TT)$. 
In order to obtain a deformation of $\Phi$ in $\wtil\M_\B(X)$, we introduce a formal power series in $t$ : 
$$
a(t)=a_1t+\frac1{2!}a_2t^2+\frac1{3!}a_3t^2+\cdots ,
\tag3-3-2$$
each $a_i$ is a section of $\CL^2(\TT)$. 
We define a formal power series $g(t)$ by 
$$
g(t)=\exp(a(t)) \in \CSpin_0(\TT)[[t]]. 
\tag3-3-3$$
The group $\CSpin_0(\TT)$ acts on differential forms and we have 
$$\align
e^{a(t)}\cdot\Phi=&\Phi+a(t)\cdot\Phi+\frac1{2!}a(t)\cdot a(t)\cdot\Phi+\cdots,\tag3-3-4\\
=&\Phi+(a_1\cdot\Phi )t+\frac1{2!}\((a_2+a_1\cdot a_1)\cdot\Phi\)t^2+\cdots.
\endalign
$$
The equation that we want to solve is, 
$$
d e^{a(t)}\cdot\Phi=0.
\tag eq$_*$ $$
At first we take $a_1$ such that $da_1\cdot\Phi=0$ as an initial condition. 
It follows from lemma 2-7 that we have 
$$
e^{-a(t)}\cdot d\cdot e^{a(t)} =\(\,(\exp(\Ad_{a(t)}) d\,\),
\tag3-3-5$$
where $(\,\exp(\Ad_{a(t)})\,)$ is the operator acting on differential forms which defined by 
the power series in $t$ : 
$$\align
\(\,(\exp(\Ad_{a(t)}) d\,\)=&d+\frac1{k!}\sum_{k=1}^\infty\Ad^k_{a(t)} d,\\
=&d+[d, a(t)]+\frac1{2!}[[d,a(t)],a(t)]+\cdots,\\
=&d+[d,a_1] t+\frac1{2!}\([d,a_2]+[[d,a_1], a_1]\)t^2+\cdots,
\endalign
$$
where $\Ad_{a(t)}^kd=[\Ad_{a(t)}^{k-1}d,\, a(t)]$.
Hence the (eq$_*$) is equivalent to the equation 
$$
\(\,(\exp(\Ad_{a(t)}) d\,\)\Phi=0,
\tag $\wtil{\text{\rm eq}}_*$ 
$$
Then it follows from proposition 2-5 that 
$\(\,(\exp(\Ad_{a(t)}) d\,\)$ is a Clifford-Lie operator of order $3$ and we have 
$$
\(\,(\exp(\Ad_{a(t)}) d\,\)\Phi\in \bE^2(X).
\tag3-3-6$$ 
From (3-3-5), we have 
$$
d e^{a(t)}\cdot\Phi=e^{a(t)}\cdot \(\,(\exp(\Ad_{a(t)}) d\,\)\Phi.
\tag3-3-7$$
We denote by $(P(t))_{[i]}$ the $i$th homogeneous part of 
a power series $P(t)$ in $t$.
Then from (3-3-7), we have 
$$
\align 
(de^{a(t)}\cdot\Phi)_{[k]}=&\sum\Sb k=i+j, \\i,j\geq 0\endSb (e^{a(t)})_{[i]}\(\,(\exp(\Ad_{a(t)}) d\,\)_{[j]}\Phi
\tag3-3-8\endalign
$$
Since $da_1\cdot\Phi=0$, we have 
$$
\(\,(\exp(\Ad_{a(t)}) d\,\)_{[0]}\cdot \Phi=\(\,(\exp(\Ad_{a(t)}) d\,\)_{[1]}\cdot \Phi=0.
\tag3-3-9$$
Thus it suffices to determine $a_k$ satisfying 
$(\wtil{\text{\rm eq}}_*)$ by induction $k$.
We assume that $a_1,\cdots, a_{k-1}$ have been determined so that 
$$
\(\,(\exp(\Ad_{a(t)}) d\,\)_{[l]}\Phi=0, \quad ( l=0, 1,\cdots, k-1).
\tag3-3-10$$
Then it follows from (3-3-8) that 
$$
(de^{a(t)}\cdot\Phi)_{[k]}=\(\,(\exp(\Ad_{a(t)}) d\,\)_{[k]}\Phi
\tag3-3-11$$
Then form (3-3-6) we see that 
$$
(de^{a(t)}\cdot\Phi)_{[k]}\in \bE^2(X).
\tag3-3-12$$
The $k$th part $(de^{a(t)}\cdot\Phi)_{[k]}$ is written as 
$$
(de^{a(t)}\cdot\Phi)_{[k]}=\frac1{k!}da_k\cdot\Phi+\Ob_k,
\tag3-3-13$$
where $\Ob_k(=\Ob_k(a_1,\cdots ,a_{k-1}))$ is the non-linear term depending only on $a_1,\cdots ,a_{k-1}$.
Since $da_k\cdot\Phi \in d\bE^1(X)\subset \bE^2(X)$, it follows from (3-3-12) that 
$$
\Ob_k\in \bE^2(X).
\tag3-3-14$$
Since $\Ob_k$ is $d$-exact, we have the cohomology class 
$[\Ob_k]\in H^2(\#_{\B})$.
Then we have 
\proclaim{Lemma 3-12} There exists a section $a_k$ satisfying 
$(de^{a(t)}\cdot\Phi)_{[k]}=0$ if and only if the class $[\Ob_k]\in H^2(\#_\B)$ vanishes. .
\endproclaim
\demo{proof} 
The equation $(de^{a(t)}\cdot\Phi)_{[k]}=0$ is written as 
$$
\frac1{k!}a_k\cdot\Phi=-\Ob_k,
\tag3-3-15$$
where $\Ob_k$ only depends on $a_1, \cdots, a_{k-1}$. 
The L.H.S of (3-3-15) is an element of the image $d\bE^1(X)$ in the complex $\#_\B$ : 
$$
\minCDarrowwidth{0.3cm}
\CD
\cdots @>d_{-1}>>\bE^0@>d_0>>\bE^1@>d_1>>\bE^2@>d_2>>\cdots.
\endCD
$$
The R.H.S. of (3-3-15) is an $d_2$-closed element of $\bE^2$ which yields the class $[\Ob_k]\in H^2(\#_\B)$. If we have $a_k$ satisfying the equation (3-3-15),
then the class $[\Ob_k]$ vanishes. 
The complex $\#_\B$ is an elliptic complex and we have the Green operator $G_{\#_\B}$ 
of the complex $\#_\B$. 
If the class $[\Ob_k]$ vanishes, we can obtain $a_k$  by using the Green operator : 
$$
\frac1{k!}a_k\cdot\Phi= -d^*G_{\#_\B}(\Ob_k)\in \bE^1.
$$ 
Then $a_k$ satisfies the equation (3-3-15).
\enddemo
We call $[\Ob_k]$ the $k$-th obstruction class. 
(Note that $[\Ob_k]$ can be defined if the lower obstruction classes vanish.) 
Since $\Ob_k$ is $d$-exact, we have that the class $[\Ob_k]\in H^2(\#_\B)$ is in the kernel of the map $p^2_\B$. 
Hence if the map $p^2_\B\: H^2(\#_\B)\to H^*_{dR}(X)$ is injective then $[\Ob_k]$ vanishes. Hence from 3-3-11, we have $a_k$ satisfying $\(\,(\exp(\Ad_{a(t)}) d\,\)_{[k]}=0$. 
By induction, we have a formal power series $a(t)$ which is a solution of the equation 
$\wtil{\text{\rm eq}}_*$. 
The rest is to show the convergence of the power series $a(t)$. The convergence can be shown essentially by the same method as in [Go]. We also have the smoothness of solutions by the standard elliptic regularity method.  Hence the result follows.
\qed\enddemo 
\head \S4. Generalized Calabi-Yau (metrical) structures \endhead 
\subhead \S 4-1 Generalized SL$_n(\C)$ structures \endsubhead 
Let $V$ be the real vector space of dim $2n$ and $\J(V)$ the set of complex structures on $V$. 
We denote by $\w^{n,0}_J V^*_\C$ the space of complex forms of type $(n,0)$ with respect to 
$J\in \J(V)$. 
Let $\P(V)$ be the set of pairs consisting of  complex structures $J$ and a non-zero complex form of type $(n,0)$: 
$$
\P(V) :=\{\, (J, \Ome_J)\, |\, J\in \J(V), \,\,0\neq \Ome_J\in \w^{n,0}_JV^*_\C\, \}.
$$
Then we have the projection to the second component  
$$\pi_2\: \P(V)\to \w^nV^*_\C.$$
\proclaim{Definition 4-1-1} 
A complex $n$-form $\Ome_V$ on $V$ is an $\SL_n(\C)$ structure if $\Ome_V$ is in the image of 
$\pi_2(\P(V))$. The set of $\SL_n(\C)$ structures on $V$ is denoted by $\A_{\SL}(V)$.
\endproclaim
Hence each SL$_n(\C)$ structure $\Ome_V$ is a complex form of type $(n,0)$ with respect to a complex  structure $J\in \J(V)$. 
Conversely for each SL$_n(\C)$ structure $\Ome_V$ 
we define a complex subspace $\ker\Ome_V$ 
by 
$$
\ker\Ome_V:=\{ v\in V_\C\, |\, i_v\Ome_V=0\, \}.
$$
Then the complexified vector space $V_\C$ is decomposed into 
$\ker\Ome_V$ and the conjugate space $ \ol{\ker\Ome_V}$ :
$$
V_\C=\ker\Ome_V\oplus \ol{\ker\Ome_V}. 
\tag 4-1-1$$
Hence we define a complex structure $J$ on $V$ by using decomposition (4-1-1) such that 
$\Ome_V$ is the complex form of type $(n,0)$ with respect to $J$. 
Then we have the map from the set of SL$_n(\C)$ structures to the set of complex structures : 
$$
\A_{\SL}(V) \to \J(V).
$$
By taking a suitable basis $\{\theta^1,\cdots,\theta^n\} $ of $\ker\Ome_V$, we can write $\Ome_V=
\theta^1\w\cdots\w\theta^n$. Then it follows that the real linear group GL$(V)$ acts on $\A_{\SL}(V)$ transitively 
with isotropy group SL$_n(\C)$ and $\A_{\SL}(V)$ is the orbit which is described as the homogeneous space:
$$
\A_{\SL}(V)=\GL(V)/\SL_n(\C).
$$
The real  conformal pin group $\CPin(\VV)$ acts on $\w^*V^*\otimes\C$. When we consider complex forms 
as pairs of real forms, we can apply the construction in section 3.
\proclaim{Definition 4-1-2}
Let $\B_{\SL}(V)$ be the orbit of $\CPin$ including SL$_n(\C)$ structures $\A_{\SL}(V)$. 
An element  $\phi_V$ of $\B_{\SL}(V)$ is a generalized SL$_n(\C)$ structure on $V$ and 
we call $\B_{\SL}(V)$ the orbit of generalized SL$_n(\C)$ structures.
\endproclaim
Let $X$ be a compact and oriented real manifold of dim $2n$. Then by applying the construction as in 
section 3, we define  $\B_{\SL}(V)$-structures on $X$ which are generalized geometric structures corresponding to the orbit $\B_{\SL}(V)$. 
Assume that there exists a $\B_{\SL}(V)$-structure $\phi$ on $X$.
Then we have the sequence of vector bundles $\{\bE^k_{\SL}\}$ over $X$ and the complex $\#_{\B_{\SL}}$ : 
$$
\minCDarrowwidth{0.3cm}
0@>>>\bE^{-1}_{\SL}@>>>\bE^0_{\SL}@>>>\bE^1_{\SL}@>>>\bE^2_{\SL}@>>>\cdots.
\tag$\#_{\B_{\SL}}$ 
$$
Let $L_\phi$ be the vector bundle over $X$ which is defined by 
$$
L_\phi=\{\, E\in \TT\, |\, \, E\cdot \phi=0\,\,\}.
$$
Then we have a decomposition : 
$$
(\TT)\otimes \C=L_\phi\oplus\ol{L_\phi},
\tag4-1-2$$
where $\ol{L_\phi}$ is the conjugate bundle of $L_\phi$.
We denote by $\w^i\ol{L_\phi}$ the $i$-th wedge product of $\ol{L_\phi}$ which acts on $\phi$ by 
the Clifford multiplication.
Then we define a vector bundle $U^i_\phi$ by 
$$
U^{-n+i}_\phi:=\w^i \ol{L_\phi}\cdot\phi,
$$
for $i=0,\cdots 2n$. The bundle $U^{-n}_\phi$ is the line bundle generated by $\phi$. 
The vector bundle $\bE_{\SL}^k$ is described in terms of $U_\phi^i$.
\proclaim{Lemma 4-1-3} 
We have the following identification as real vector bundle :
$$\align 
&\bE_{\SL}^0\cong U_\phi^{-n+1},\\
&\bE_{\SL}^1\cong U_\phi^{-n}\oplus U_\phi^{-n+2},\\
&\bE_{\SL}^2\cong  U_\phi^{-n+1}\oplus U_\phi^{-n+3}.
\endalign
$$
In general we have 
$$\align
&\bE_{\SL}^{2k-1}\cong  \oplus_{i=0}^k U_\phi^{-n+2i},\\
&\bE_{\SL}^{2k}\cong \oplus_{i=0}^kU_\phi^{-n+2i+1}.
\endalign
$$
\endproclaim
\demo{proof} We consider the complex form $\phi=\phi^{\Re}+\sqrt{-1}\phi^{\Im}$ as  the pair of real forms 
$(\phi^{\Re}, \phi^{\Im})$. Then applying the construction in section 3, we have the vector bundles $\bE^k_{\SL}$ which generated by  
$$
\bE^k_{\SL}=\text{\rm span}\{\, (a\cdot\phi^{\Re}, a\cdot\phi^{\Im})\, |\, a\in \CL^k\,\}.
$$
Then we have the complex form $a\cdot\phi^{\Re}+\sqrt{-1}a\cdot\phi^{\Im}=a\cdot\phi$. 
From the decomposition (4-1-2), we have the identification : 
$$
\CL^{2k}\otimes\C\cong \CL^{2k}(L_\phi\oplus\ol{L_\phi})\cong\oplus_{i=0}^k\w^{2l}(L_\phi\oplus\ol{L_\phi}).
$$
Since $L_\phi\cdot\phi=\{0\}$,  We have an identification : 
$$\align
\bE^{2k-1}_{\SL}=\CL^{2k}\cdot\phi\cong&\oplus_{l=0}^k\w^{2l}\ol{L_\phi}\cdot\phi\\
=&\oplus_{l=0}^kU_\phi^{-n+2l}.
\endalign
$$
Similarly we have $\bE_{\SL}^{2k}\cong  \oplus_{i=0}^kU_\phi^{-n+2i+1}$.
\qed\enddemo
\proclaim{Propositon 4-1-4} 
The complex $\#_{\B_{\SL}}$ is elliptic, so that is, the orbit $\B_{\SL}$ is an elliptic orbit.
\endproclaim
\demo{proof} Since there is the inclusion $\CL^{k-2}\subset\CL^k$ , we have the inclusion $\bE_{\SL}^{k-2}\subset\bE_{\SL}^k$ with the quotient 
$$
\bE^k_{\SL}/\bE^{k-2}_{\SL}\cong U_\phi^{-n+k+1}, 
$$
for $k\geq 0$. 
Replacing $\bE^{-1}$ by $\bE^{-1}\otimes\C$, we have a complex $\wtil\#_{\SL}^*$
Hence there is a map of the complex $\wtil\#_{\B_{\SL}}$ by shifting its degree from $*$ to $*+2$ : 
$$
\wtil\#_{\B_\SL}\overset[2]\to\longrightarrow\wtil\#_{\B_{\SL}}.
$$
Hence we have the following commutative diagram : 
$$\minCDarrowwidth{0.3cm}
\CD
 @.0@>>>0@>>>\bE^{-1}_{\SL}\otimes\C @>>> \bE^0_{\SL}@>>>\bE^1_{\SL}@>>>\cdots \\
@. @VVV @VVV @VVV @VVV @VVV\\
 0@>>>\bE^{-1}_{\SL}\otimes\C @>>> \bE^0_{\SL}@>>>\bE^1_{\SL} @>>> \bE^2_{\SL}@>>>\bE^3_{\SL}@>>>\cdots \\
  @. @VVV @VVV @VVV @VVV@VVV\\
0@>>> U^{-n}_\phi @>\ol\pa>> U_\phi^{-n+1}@>\ol\pa>> U_\phi^{-n+2}@>\ol\pa>>U_{\phi}^{-n+3}@>\ol\pa>>U_\phi^{-n+4}@>>>\cdots.
\endCD
$$
It follows from $U_\phi^{-n+i}=\w^i\ol{L_\phi}\cdot\phi$ that the quotient complex $(U^{-n+*}_\phi,\ol\pa)$ is an elliptic complex. Hence the from the commutative diagram we see that the complex $\#_{\B_{\SL}}$ is elliptic by induction on degree $k$.
\qed\enddemo
The complex $(U_\phi^p,\ol\pa)$ is the deformation complex of generalized complex structures which is introduced by [Gu1]. 
Then the exterior derivative $d$ acting on $U_\phi^p$ is decomposed into two projections $\pa$ and $\ol\pa$, so that is, 
$$d=\pa+\ol\pa,
$$ 
$$
U_\phi^{p-1}\overset\pa\to{\longleftarrow }U_\phi^p\overset\ol\pa\to{\longrightarrow} U_\phi^{p+1}.
$$
Let $\Phi$ be a $\B_{\SL}(V)$- structure on $X$. Then there is the generalized complex structure $\J_\phi$ corresponding to $\phi$. 
We define an operator $d^\J$ by 
$$
d^\J:= \sqrt{-1}(\ol\pa-\pa).
$$
The $dd^\J$-property is introduced by [Hi], [Gu2] and [Cav]: 
\proclaim{Definition 4-1-5} 
A generalized complex manifold $(X, \J)$ satisfies the $dd^\J$-property iff the following are equivalent : 
\roster 
\item $a\in \w^*T^*$ is $d$-closed and $d^\J$-exact,
\item $a\in \w^*T^*$ is $d$-exact and $d^\J$-closed, 
\item $a=dd^\J b$ for some $b\in\w^*T^*$.
\endroster 
\endproclaim
\proclaim{Theorem 4-1-6} 
If the $dd^\J$-property holds for the $\J_\phi$ corresponding a $\B_{\SL}$-structure $\phi$, then 
$\phi$ is a topological structure, so that is, we have unobstructed deformations of $\phi$
on which the local Torelli type theorem holds.
\endproclaim
\demo{proof} Since $U_\phi^p$ is the eigenspace of the action of $\J_\phi$ with eigenvalue 
$\sqrt{-1} p$. Hence 
we have the decomposition $\w^*T^* =\oplus_{p=-n}^n U_\phi^p$.
If an exact form $da^{(m)}$ is an element of $U_\phi^{m-1}$ for $a^{(m)}\in U_\phi^m$, we have 
$\pa da^{(m)} =\pa\ol\pa a^{(m)}=0$. 
Hence applying the $dd^\J$-property we have 
$$\align
da^{(m)}=&dd^\J b =2\sqrt{-1}\,\ol\pa\pa b=2\sqrt{-1}d\pa b,
\endalign
$$ for 
$b\in U_\phi^{m-1}$. 
Then we have $da^{(m)} =d\gam$ for $\gam =2\sqrt{-1}\pa b\in U_\phi^{m-2}$.
From our decomposition, a form $a$ is written as 
$$
a=\sum_{p=-n}^m a^{(p)},
$$
where $a^{(p)}\in U_\phi^p$ for some $m$. 
If $da$ is an element of $\sum_{p=-n}^kU_\phi^p$, then applying the $dd^\J$-property successively, 
we have $da=db$ for $b \in \sum_{p=-n}^{k-1} U_\phi^p$. 
Similarly if $da\in \w^{\even}T^*$ ( resp. $da\in \w^{odd}T^*$) then
applying the property we see that 
$da=db$ for $b\in \w^{odd}T^*$ (resp. $b\in \w^{\even}T^*$).
Hence it follows from lemma 4-1-3 that if $da\in \bE_{\SL}^k$ then $da=db$ for 
$b \in \bE_{\SL}^{k-1}$ ($k\geq 1$). 
It implies that the map $p^k_{\B}\: H^k(\#_{\B_{\SL}})\to H^*_{dR}(X)$ is injective for 
$k\geq 1$. 
\qed\enddemo
Gualtieri also shows that the $dd^c$-property holds for generalized K\"ahler structures [Gu2]. 
By applying his theorem, we have
\proclaim{Theorem 4-1-7} 
Let $\phi$ be a $\B_{\SL}(V)$-structure on $X$ with the corresponding generalized complex structure $\J_\phi$. 
if there exists another generalized complex structure $\I$ such that 
the pair $(\I, \J_\phi)$ defines a generalized K\"ahler structure on $X$, then 
$\B_{\SL}(V)$-structure $\phi$ is a topological structure.
\endproclaim
As in the proof of proposition 4-1-4, we have the short exact sequence : 
$$
\CD 
0@>>> \wtil\#_{\B_{\SL}}@>[2]>> \wtil\#_{\B_{\SL}}@>>> (U_\phi^*,\olpa)@>>>0.
\endCD
$$
It follows from $\bE^0_{\SL}\cong U_\phi^{-n+1}$ that we have the long exact sequence : 
$$\minCDarrowwidth{0.1cm} 
\CD
0@>>>H^{-1}(\wtil\#_{\B_{\SL}})@>>>H^1(\wtil\#_{\B_{\SL}})@>>>H^2_{\ol\pa}(U_\phi^*)
@>>>H^0(\wtil\#_{\B_{\SL}})@>>>H^2(\#_{\B_{\SL}}),
\endCD
$$
where $H^2_{\ol\pa}(U_\phi^*)$ denotes the cohomology group,
$$ 
H^2_{\ol\pa}(U_\phi^*)=\(\ker\ol\pa \: U_\phi^{-n+2}\to U_\phi^{-n+3} \)/\ol\pa (U_\phi^{-n+1}).
$$ 
which is the infinitesimal tangent space of deformations of generalized complex structures 
( see Chapter 5 in [Gu1]). 
Since the complex $\wtil\#_{\B_{\SL}}$ is a subcomplex of the complexified de Rham complex we have the map $\til p^k_\B\: H^k(\#_{\B_{\SL}})\to H^*_{dR}(X,\C)$ as in section 3. Then we have a diagram : 
$$
\CD
H^0(\wtil\#_{\B_{\SL}})@>>> H^2(\wtil\#_{\B_{\SL}})\\
@V \til p^0_\B VV  @VV \til p^2_\B V \\
H^*_{dR}(X) @>\cong >> H^*_{dR}(X).
\endCD
$$
Hence it implies that if $\til  p^0_\B$ is injective then the map 
$H^0(\wtil\#_{\B_{\SL}})\to H^2(\wtil\#_{\B_{\SL}})$ is injective. 
As in proof of theorem 4-1-6, the $dd^\J$-property implies the injectivity of the map $\til p^0_\B$. Hence we have 
\proclaim{Proposition 4-1-8}
Let $\phi$ be a generalized SL$_n(\C)$ structure on $X$. If the generalized complex structure $\J_\phi$ satisfies the $dd^\J$-property then we have the exact sequence : 
$$
0@>>> H^{-1}(\wtil\#_{\B_{\SL}})@>>>H^1(\wtil\#_{\B_{\SL}})@>>> H^2_{\ol\pa}(U_\phi^*)@>>>0.
$$
Hence the infinitesimal tangent $H^2_{\ol\pa}(U_\phi^*)$ gives rise to small deformations of 
generalized complex structure $\J_\phi$ which correspond to deformations of generalized SL$_n(\C)$ structures.
\endproclaim
\subhead \S4-2. Generalized Calabi-Yau (metrical) structures\endsubhead 
Let $\Ome_V$ be an SL$_n(\C)$-structure and $\ome_V$ a real $2$-form on the real vector space of $2n$ dim. As in section 4-1, the SL$_n(\C)$-structure $\Ome_V$ gives rise to the complex structure $J$ on $V$ and we define a bilinear form $g$ by 
$$
g(u,v) = \ome(Ju, v), \qquad ( u,v \in V)
$$
\proclaim{Definition 4-2-1} 
A pair $(\Ome_V,\ome_V)$ is a Calabi-Yau structure on $V$ if the followings hold :
\roster
\item $$\Ome_V\w\ome_V=0,$$
\item $$\Ome\w\ol\Ome= c_n\ome^n$$, 
\item The corresponding bi-linear form $g$ is positive-definite. 
\endroster
\endproclaim
The condition (1) implies that $\ome_V$ is a form of type $(1,1)$ with respect to $J$ and 
then it follows from (3) that $\ome_V$ is a Hermitian form. 
The equation (2) is called the Monge-Amp$\grave{\text{\rm e}}$re condition. 
Let $\A_{\CY}(V)$ be the set of Calabi-Yau structures on $V$ which consist of complex $n$-forms and real $2$-forms. Then 
the real linear group $GL(V)$ acts on $\A_{\CY}(V)$ transitively with the isotropy group 
SU$(n)$. 
Hence $\A_{\CY}(V)$ is the orbit of $\GL(V)$ which is described as a homogeneous space : 
$$
\A_{\CY}(V)=\GL(V)/\SU(n).
$$Let $(\Ome_V, \ome_V)$ be a Calabi-Yau structure on $V$. Then we consider a pair 
$(\Ome_V, e^{\sqrt{-1}\ome_V})$ consisting two generalized SL$_n(\C)$ structures 
$\Ome_V$ and $e^{\sqrt{-1}\ome_V}$. 
\proclaim{Definition 4-2-2} 
The orbit $\B_{\CY}(V)$ of $\CPin(\VV)$ through the pair $(\Ome_V, e^{\sqrt{-1}\ome_V})$
is called the generalized Calabi-Yau orbit. 
An element $(\phi_{\ss V,0}, \psi_{\ss V,0})$ of the orbit $\B_{\CY}(V)$ is a generalized Calabi-Yau structure on $V$. 
Note that the orbit $\B_{\CY}(V)$ is embedded into pairs of complex forms $\w^*V^*_\C\oplus \w^*V^*_\C$.
Let $X$ be a compact real manifolds of dim $2n$. 
Then as in section $3$, we define generalize Calabi-Yau (metrical) structures on $X$ as $\B_{\CY}(V)$-structures on $X$.
\endproclaim
Let $(\phi_0,\phi_1)$ be a generalized Calabi-Yau structure on $X$. Since it consists on 
generalized SL$_n(\C)$ structures, we obtain the pair $(\J_0, \J_1)$ of the corresponding generalized complex structures on $X$. 
Then we see that the pair $(\J_0,\J_1)$ is a generalized K\"ahler structure. 
By applying Gualtieri's theorem, we obtain the following theorem of deformations of generalized 
Calabi-Yau structures (which are deformations of pairs consisting two generalized $\SL_n(\C)$ structures with the conditions): 
\proclaim {Theorem 4-2-3} 
The generalized Calabi-Yau orbit is an elliptic and topological orbit.
\endproclaim
\demo{proof} Let $(\phi_0,\phi_1)$ be a generalized Calabi-Yau (metrical) structure with 
the generalized K\"ahler structure $(\J_0,\J_1)$ on $X$. 
We denote by $\#_{\B_{\CY}}=\{\bE_{\CY}^*, d\}$ the deformation complex of generalized Calabi-Yau structure 
$(\phi_0, \phi_1)$. 
Then it suffices to show that each map 
$$p^k_{\B_{\CY}}\: H^k(\#_{\B_{\CY}})\to \oplus^2H^*_{dR}(X,\C)
$$ 
is injective for $k=1,2$.
We have the eigenspace decomposition of $\w^*T^*=\oplus_{p=-n}^n U_{\phi_i}^p$ for each $i=0,1$.
Since $[\J_0,\J_1]=0$, we have a further decomposition : 
$$
\w^*T^*=\bigoplus\Sb |p+q|\leq n \\ p+q\equiv n (\text{\rm mod} 2)\endSb U^{p,q},
$$
where $U^{p,q}=U^p_{\phi_0}\cap U^q_{\phi_1}$. 

Each $\bE^k_{\CY}$ consists of pairs of complex forms. 
Then the projection $\pi_1$ to the first component induces a map from the complex 
$\#_{\B_\CY}$ to $\#_{\B_\SL}$. We denote by $(\ker^*,d)$ the complex defined by the kernel of $\pi_1$. Then we have a short exact sequence : 
$$\minCDarrowwidth{0.3cm}
\CD
0@>>>(\ker^*,d)@>>>\#_{\B_\CY}@>>>\#_{\B_\SL}@>>>0,
\endCD
\tag4-2-1$$
so that is, 
$$\minCDarrowwidth{0.3cm}
\CD
@.\ker^0@>>>\ker^1@>>>\ker^2@>>>\cdots\\
@. @VVV @VVV @VVV \\
\bE_\CY^{-1}@>>> \bE_{\CY}^0@>>> \bE_{\CY}^1@>>>\bE_{\CY}^2@>>>\cdots\\
@VVV @VVV @VVV @VVV \\
\bE_\SL^{-1}@>>>\bE_{\SL}^0@>>>\bE_{\SL}^1@>>>\bE_{\SL}^2@>>>\cdots.\\
\endCD
$$
If $E\cdot\phi_1=0$ for real $E\in \TT$ then we see that $E=0$. It implies that 
$\ker^{0}\cong\{0\}$. Similarly $\ker^1$ and $\ker^2$ are respectively given by 
$$\align
&\ker^1\cong U^{0,-n+2},\\
&\ker^2\cong U^{1,-n+1}\oplus U^{-1,-n+1}\oplus U^{1,-n+3}\oplus U^{-1,-n+3}.
\endalign
$$
The complex $(\ker^*, d)$ is a subcomplex and we have the map $p^k_{\ker}\: 
H^k(\ker^*)\to H^*_{dR}(X)$. 
By applying the Hodge decomposition of generalized K\"ahler manifold in [Gu2], 
we see that $p^k_{\ker}$ is injective for each $k$. 
The short exact sequence  (4-2-1) is a subsequence of the following short exact sequence 
define by the de Rham complex (dR) which yields a splitting long exact sequence, 
$$\minCDarrowwidth{0.3cm}
\CD 
0@>>> (dR)@>>>(dR)\oplus (dR)@>>>(dR)@>>>0.
\endCD
$$
Hence we have the diagram of long exact sequences:
$$
\minCDarrowwidth{0.3cm} 
\CD 
\cdots @>>>H^k(\ker^*)@>>>H^k(\#_{\B_\CY})@>>>H^k(\#_{\B_\SL})@>>>\cdots\\
@. @Vp^k_{\ker}VV @Vp^k_{\B_\CY}VV @Vp^k_{\B_{\SL}}VV @. \\
0@>>> H^*_{dR}(X)@>>>\oplus^2 H^*_{dR}(X)@>>>H^*_{dR}(X)@>>>0,
\endCD
$$
where the sequence at top is the long exact sequence of (4-2-1). 
Since $p^k_{\B_\SL}$ is injective, we see that $p^k_{\B_{\CY}}$ is also injective for $k$. 
Hence the results follows.
\qed\enddemo
\head \S 5. Genealized hyperK\"ahler, G$_2$ and Spin$(7)$ structures \endhead
\subhead \S5-1. Generalized hyperK\"ahler structures
\endsubhead
In this section we will introduce two types of generalized hyperK\"ahler structures, so that is, 
type 1 and type 2.
A genralized hyperK\"ahler structure of type 1 is defined by closed differential forms as in section 3. The one of type 2 is based on generalized complex structures satisfying a relations.
Let $V$ be a $4m$ dimensional real vector space. A hyperK\"ahler structure on $V$ with 
a hyperK\"ahler structure $(g, I,J,K)$, where $I,J$ and $K$ are three complex structure satisfying the quoternion relations, ($I^2=J^2=K^2=IJK=-1$) and $g$ is a Hermitian metric with respect to $I,J$ and$K$. Then we have three Hermitian form $\ome_I, \ome_J$ and $\ome_K$ with respect to $I, J$ and $K$. Conversely such three forms $(\ome_I, \ome_j,\ome_K)$ yields the hyperK\"ahler structure $(g, I,J,K)$. Hence hyperK\"ahler structures on $V$ can be regarded as geometric structures defined by three special $2$-forms (see [Go]). 
As in section $4$, three generalized SL$_n(\C)$ structures are defined by  
$$
\phi_I:=e^{\sqrt{-1}\ome_I}, \quad\phi_J:=e^{\sqrt{-1}\ome_J}, \quad\phi_K:=e^{\sqrt{-1}\ome_K}.
$$
\proclaim{Definition 5-1-1} Let $\B_{\HK}(V)$
be the orbit of $\CPin$ through the triple $(\phi_I, \phi_J, \phi_K)$. 
We call $\B_{\HK}(V)$ the orbit of generaliazed hyperK\"ahler structures on $V$. 
A generalized hypK\"ahler structure of type 1 on a $4m$-fold $X$ is a $B_{\HK}(V)$-structure on $X$.
\endproclaim
It is worthwhile to mention that there exist six generalized complex structures defined by 
generalized hypK\"ahler structures. 
As in section $4$, generalized SL$_n(\C)$ structures yield generalized complex structures. 
Hence we have three generalized complex structures $\I_1$, $\J_1$ and $\K_1$ corresponding to $\ome_I, \ome_J $ and $\ome_K$ respectively. 
It immiditately follows that three compositions $\K_0:=\I_1\J_1$, 
$\I_0:=\J_1\K_1$ and $\J_0:=\K_1\I_1$ give generalized complex structures respectively. Hence we have six generalized complex structures 
$\I_0, \J_0, \K_0$, $ \I_1, \J_1, \K_1$ which satisfying relations : 
$$\align
&\I_0^2=\J_0^2=\K_0^2=\I_0\J_0\K_0=-1, \tag5-1-1\\
&\I_0\I_1=\I_1\I_0=\J_0\J_1=\J_1\J_0=\K_0\K_1=\K_1\K_0=-G,\tag5-1-2\\
\endalign
$$
where $G$ is a section of the bundle SO$(\TT)$ with $G^2=1$, which is called a generalized metric. 
The relation (5-1-1) implies that the triple $(\I_0, \J_0,\K_0)$ satisfies the quaternion relations and (5-1-2) shows that we have three generalized K\"ahler structures with a same $G$. 
\proclaim{Definition 5-1-2} 
A generalized hyperK\"ahler structure of type 2 is a system consisting of six generalized complex structures with the relation (5-1-1,1).
\endproclaim
The algebra generated by $\lan G, 1\ran$ with $G^2=1$ is isomorphic to the direct sum 
$\R\oplus \R$ by taking a basis,
$$
\frac1{\sqrt2}(1+G), \quad \frac1{\sqrt2}(1-G).
$$
When we set 
$$\I_1=-G\otimes i,\quad  \J_1=-G\otimes j,\quad \K_1=-G\otimes k,
$$
then we see that the algebra generated by $\I_\a,\J_\a, \K_\a, G, 1, \a=0,1$ with relations 
(5-1-1,2) is isomorphic to $(\R\oplus\R)\otimes\Bbb H\cong \Bbb H\oplus \Bbb H$. 
In the case of generalized K\"ahler structure, we have the algebra $(\R\oplus\R)\otimes\C\cong \C\oplus\C$.
Hence there are following correspondence : 
$$\matrix
\text{\rm generalized} &\text{\rm generalized }
&\text{\rm generalized}\\
\text{\rm  metric } &\text{ \rm K\"ahler structure }&\text{\rm  hyperK\"ahler structure} \\
\updownarrow &\updownarrow&\updownarrow\\
\R\oplus\R&\C\oplus \C&\Bbb H\oplus\Bbb H
\endmatrix
$$
\subhead \S 5-2. Generalized G$_2$ structures
\endsubhead
Let $\O$ be the octernions which is regarded as the $8$ dimensional real vector space with the metric
$\lan\, ,\,\ran$  and 
the non associative multiplication. Consider the $7$ dimensional vector space $W:=$Im$\O$ with a volume form vol$_7$ and define a $3$ form $\phi_{\G}$ on $W$ by 
$$
\phi_{\G}(x,y,z) := \lan xy, z\ran.
$$
A $4$ form $\psi $ on $W$ is defined by the Hodge star operator $\star$, 
$$
\psi_{\G}:=\star\phi_{\G}.
$$
Let $\B_{\G}$  be the orbit of $\CPin$ through 
the pair ($\vol_7-\phi_{\G}, \, 1-\psi_{\G}$). Then we have
\proclaim{Definition 5-2} 
A generalized $\G$ structure is a $\B_{\G}$-structure on a $7$-fold $X$.
\endproclaim
\subhead 
\S5-3. Generalized Spin$(7)$ structures
\endsubhead
As in section 5-2, we consider $\O$ as the $8$ dimensional vector space $V$ which is decomposed into $V=\R\oplus$Im$\O$, where $\R$ denotes the real part of the octernions 
with a non-zero one form $e^1$. Then the Spin(7) form $\phi_{\Spin}$ (the Cayley
form) is a $4$-form on $V$ defined by 
$$
\phi_{\Spin}=e^1\w\phi_{\G}+\psi_{\G}.
$$
Let $\B_{\Spin}$ be the orbit of the group $\CPin(\VV)$ through
$$\Phi_{\Spin}:=1-\phi_{\Spin}+\vol_8,$$ where $\vol_8$ denotes the volume form of $V\cong\O$.
Then we have
\proclaim{Definition 5-3} 
A generalized $\Spin(7)$ structure is a $\B_{\Spin}$-structure on a $8$-manifold $X$.
\endproclaim
Let $(\Ome_I, e^{\sqrt{-1}\ome_I})$ be a Calabi-Yau structure on a $8$-manifold $X$. Then 
the real part $\Phi=\Ome_I^{\Re}+(e^{\sqrt{-1}\ome_I})^{\Re}$ yields a generalized 
Spin$(7)$ structure on $X$. Since the group $\CPin$ is real,  for a generalized Calabi-Yau structure $(\phi_0, \phi_1)$  on a $8$-manifold $X$, the real part 
$$
\Phi=\phi_0^{\Re}+\phi_1^{\Re}
$$
is a generalized Spin$(7)$ structure.
\proclaim{Proposition 5-4} Let $X$ be a compact, oriented $8$-manifold with a generalized Spin$(7)$ structure $\Phi$.
Deformations of generalized Spin$(7)$ structures on $X$ are unobstructed.
\endproclaim
Our proof is a generalization of the proof of deformations of Spin$(7)$ structures in [Go]. 
Let $\CSpin_0(\VV)$ be the identity component of  the conformal spin group which acts on $\w^*T^*$. There is the metric $g_{\ss\O}$ on the octernions $\O$ which yields a generalized metric $G_{\ss\O}$ by 
$$
G_{\ss\O}=\pmatrix
o&g^*_{\ss\O}\\
g_{\ss\O}&0
\endpmatrix,
$$
where $g^*$ is the dual metric on the dual space $\O^*$. 
Then we have 
\proclaim{Lemma 5-5} 
Let $H$ be the isotropy group which is defined by  
$$
H:=\{ g\in \CSpin_0(\VV)\,|\, g\cdot\Phi_{\Spin}=\Phi_{\Spin}\,\}.
$$Then the identity component  $H_0$ of $H$ preserves the generalized metric $G_{\ss\O}$.
\endproclaim
\demo{proof} 
Every $g\in \CSpin_0(\VV)$ is written as $g=e^a$ for $a\in \CL^2$. 
Then $\CL^2$ is decomposed into 
$$
\CL^2=\R\oplus \text{\rm End}(V)\oplus \w^2V^*\oplus\w^2V.
$$
Then $a$ is written as $a=\lam +A+b+\b$ for  $A\in $End$(V)$, $\b\in \w^2V$, 
$\b\in \w^2V$.
If $a\cdot \Phi_{\Spin}=0$ for $a=\lam +A+\b +b$, then we see that 
$$
\align 
&(\lam+A)\cdot\Phi_{\Spin}=0,\tag5-3-1\\
&(b+\b)\cdot\Phi_{\Spin}=0\tag5-3-2
\endalign
$$
From (5-3-1), we have 
$$\lam=0\quad \text{\rm  and }\quad A\in \text{\rm spin}(7).\tag3-5-3$$
The space of $2$-forms and the space of $2$-vectors are respectively decomposed into 
representation spaces of Spin$(7)$,
$$\align
&\w^2 V^*=\w^2_7V^*\oplus\w^2_{21}V^*,\tag5-3-4\\
&\w^2 V=\w^2_7 V\oplus\w^2_{21} V,\tag5-3-5
\endalign
$$
where $\w^2_{21}V^*\cong \w^2_{21}V$ is isomorphic to the Lie algebra spin$(7)$ and 
$\w^2_7$ denotes the irreducible $7$ dimensional representation space.
Let $p^*$ be the dual $2$-vector of $p\in \w^2_7V^*$ and $q^*$ the dual $2$-vector 
of $q\in \w^2_{21}V^*$.
Then we have 
$$\l\{
\aligned
&p\w\Phi_{\Spin}=p-3\star p,\\
&p^*\cdot\Phi_{\Spin}=3p-\star p,\\
&q\w\Phi_{\Spin}=q+\star q,\\
&q^*\cdot\Phi_{\Spin}=-(q+\star q),
\endaligned
\r\}\tag5-3-6
$$ 
where $\star$ denotes the Hodge star operator with respect to the metric $g_{\ss\O}$.
From (5-3-6) if $(b+\b)\cdot\Phi_{\Spin}=0$, then 
$$b+\b= q+q^*,\quad \(\text{\rm  for }q \in 
\w^2_{21}\cong \text{\rm spin}(7)\).\tag5-3-7$$
From (5-3-3,7) the Lie algebra of the isotropy group $H$ is given by the direct sum 
spin$(7)\oplus$spin$(7)$. Hence the identity component $H_0$ is the product Spin$(7)\times$Spin$(7)$ which preserves the generalized metric $G_\O$.
\qed\enddemo
Let $\CSpin_0(X)$ be the fibre bundle over a compact, oriented manifold $X$ with fibre 
$\CSpin_0$. It suffices to consider that a generalized Spin$(7)$ structure $\Phi$ is a closed section of the fibre bundle $\CSpin_0(X)$.
For each $x\in X$, we take an identification $h:\O\cong T_xX $ which gives the identification $\h{h}\: T_xX\oplus T^*_xX\cong \O\oplus \O^*$. Then there is a $e^a\in \CSpin_0(\VV)$ such that 
$$
\Phi_{\Spin}=e^a\cdot h^*\Phi\in \w^*\O^*.
$$
We denote by $\Ad_{e^a}$ the adjoint of $e^a$. Then the composition $\h{h}\circ\Ad_{e^a}$
gives the identification $\O\oplus \O^*\cong T_xX\oplus T^*_xX$ which defines a generalized metric $G_{\Phi_x}$ by 
$$
G_{\Phi_x}:=\h{h}_*\circ (\Ad_{e^a})_* G_\O.
$$
When we take an other identification $h'\: \O\cong T_xX$ preserving the orientation and an element $e^{a'}\in \CSpin_0(\VV)$ with 
$$
\Phi_{\Spin}=e^a\cdot h^*\Phi=e^{a'}\cdot(h')^*\Phi,
$$
we have $h^{-1} e^a\Phi_{\Spin}=(h')^{-1}e^{-a'}\Phi_{\Spin}$. 
We can take $e^{a'}h'h^{-1} e^{-a}$ is an element of identity component.
It follows from lemma 5-5 that 
$$e^{a'}h'h^{-1} e^{-a}\in H_0\cong \Spin(7)\times\Spin(7).
$$
Hence 
$$
G_{\Phi_x}:=\h{h}_*\circ (\Ad_{e^a})_* G_\O=\h {h'}_*\circ(\Ad_{e^{a'}} )_* G_\O.
$$
Thus we have  
\proclaim{Lemma 5-6} 
Let $\Phi$ be a generalized Spin$(7)$ structure on $X$. Then there is a generalized metric $G_\Phi$ on $X$ which is canonically defined by $\Phi$.
\endproclaim
A generalized metric $G_\Phi$ yields an operator $*$ acting on differential forms with $*^2=1$ and 
a generalized metric $G_\Phi$ also yields a metric $\lan\,,\,\ran$ on forms and we have the adjoint operator $d^*$ which is given by $d^*=*d*$. By using adjoint operator we have the Laplacian $\triangle_{\Phi}$ by 
$\triangle_{\Phi}=d d^*+ d^* d$. (see [Gu2] for detail).
\proclaim{Lemma 5-7} 
Let $\w^{\even}_-$ be anti-self dual forms of even type with respect to the operator $*$, so that is, 
$$
\w^{\even}_-=\{\,  s\in \w^{\even}\, |\, *s=-s\,\}.
$$
Then the bundle $\w^{\even}_-$ is the subbundle of the vector bundle $\bE^1_{\Spin}(X)$
$=\CL^2\cdot\Phi$ 
 \endproclaim
 \demo{proof}
 At first we will show the lemma for the form $\Phi_{\Spin}=1-\phi_{\Spin}+\vol_8$ on the octernions $\O\cong V$. 
 The action of $g\in \GL(V)$ on forms is given by 
 $$
 (\det g)^{\frac12}\rho_{g^{-1}},
 $$
 where $\rho$ denotes the linear representation of $\GL(V)$. 
 Hence for $\lam \in \R$, $\lam1\in \GL(V)$ acts on $\Phi_{\Spin}$ by 
 $$
 \lam^4+\Phi_{\Spin}+\lam^{-4}\vol_8.
 $$
 By differentiating to $\lam$ we have that 
 $$-1+\vol_8\in \CL^2\cdot\Phi_{\Spin}.\tag5-3-8$$
 
 Note that the the Hodge star $\star$ is given by the $*$ operator, 
 $$
 \star a=\cases 
 &+ * a, \quad  a\in \w^0+\w^1+\w^4+\w^7+\w^8,\\
 &-*a,\quad a\in \w^2+\w^3+\w^5+\w^6.
 \endcases
 \tag5-3-9$$
 From (5-3-6), we see that 
 $$\w^2_7\oplus\w^6_7\oplus(\w^2_{21}\oplus\w^6_{21})_-\subset \CL^2\cdot\Phi_{\Spin},
 \tag5-3-10$$
 where $(\w^2_{21}\oplus\w^6_{21})_-=\{\, q+\star q\, |\, q\in\w^2_{21}\,\}$.
 From (5-3-9) we have $\star q=-*q$ for $q\in\w^2V^*$.
 Hence $q+\star q=q-*q$ is the anti-self dual form with respect to $*$.
 We also have 
 $$gl(8)/spin(7)= so(8)/spin(7)\oplus \{\lam 1\,|\lam \in \R\}\oplus Sym_0(8),
 $$
 where Sym$_0(8)$ denoted the trace-free symmetric matrix which yields anti-self-dual $4$-forms, so that is, 
 $$
 \w^4_-=\{\, \hrho_a\Phi \,|\, a\in Sym_0(8)\, \},
 \tag5-3-11$$
 where $\hrho$ denotes the differential representation of End$(V)$.
 Hence $\w^4_-\in \CL^2\cdot\Phi_{\Spin}$.
 The anti-self -dual form of even type $\w^{\even}_-$ is decomposed into there parts : 
 $$
 \w^{\even}_-=(\w^{0,8})_-\oplus (\w^{2,6})_-\oplus(\w^4_-),
 $$
 where $\w^{0,8}_-$ and $\w^{2,6}_-$ respectively denotes the anti-self dual $0,8$-forms and 
 anti-self dual $2,6$-forms. 
 Hence from (5-3-8,10,11) it follow that 
 $$\w^{\even}_-\subset \CL^2\cdot\Phi_{\Spin}.\tag5-3-12$$
 Every generalized Spin$(7)$ form $\Phi_V$ on $V$ is given by $ k\cdot \Phi_{\Spin}$ for 
 an element $k \in \CSpin(\VV)$. Then the action of $k$ on forms yields the identification, 
 $$
 \bE_{\Phi_{\Spin}}^1\cong\bE_{\Phi_V}^1.
 $$
 Then the operator $*_{\ss\Phi_V}$ with respect to $\Phi_V$ is given by the adjoint action of $k$, 
 $$
 *_{\ss\Phi_V}= k\cdot( *_{\Phi_{\Spin}} )\cdot k^{-1}.
 $$
 We denote by $\w^{\even}_{\ss\Phi_V}$ the anti-self-dual forms of even type with respect to 
 $\Phi_V$.
 Then the action of $k$ also induces the identification between anti-self-dual forms 
 $\w^{\even}_{\ss-,\Phi_{\Spin}}\cong \w^{\even}_{\ss -,\Phi_V}$. Then from $\w^{\even}_{\ss-,\Phi_{\Spin}}\subset \bE^1_{\ss\Phi_{\Spin}}$ we have $\w^{\even}_{\ss-,\Phi_V}\subset \bE^1_{\ss\Phi_V}$.
 Hence we have the result.
 \qed\enddemo
\demo{proof of proposition 5-4} 
Let $\Phi$ be a generalized Spin$(7)$ structure on an oriented and compact $8$-manifold $X$ with the differential complex $\#_{\B_{\Spin}}=\{\bE^k_{\Spin}(X), d_k\}$ over $X$.
It suffices to show that the map $p^2_{\B_{\Spin}}$ is injective.
Note that $\bE^1_{\Spin}(X)\subset \w^{even}$ and $\bE^2(X)\subset \w^{\odd}$. 
We have the Lapalacian $\triangle_{\Phi}\: \w^{\even/\odd}\to\w^{\even/\odd}$ with respect to the generalized metric $G_\Phi$ as in before.
We denote by $\Cal G_\Phi$ the Green operator of the Laplacian $\triangle_\Phi$ on the manifold $X$. 
Let $d\a$ be a $d$-exact section of $\bE^2(X)$. Then by using the Hodge decomposition with respect to $\triangle_\Phi$, we have 
$$
d\a = d( dd^*\Cal G_\Phi \a + d^* d \Cal G_\Phi \a).
$$
We define $\a_-\in \w^{\even}_-$ by 
$$
\a_-=d^*d\Cal G_\Phi-*d^*d\Cal G_\Phi.
$$
Then we  have $d\a =d\a_-$. 
It follows from lemma 5-7 that $\a_-\in \bE^1_{\Spin}(X)$ and it implies that the map $p^2_{\B_{\Spin}}$ is injective.
\qed\enddemo 
\proclaim{Acknowledgements} 
 {\rm The author is grateful to Professor Nijel Hitchin for nice discussions on the generalized geometry
 and he is thankful to Professor Akira Fujiki for useful conversations and encouragements. 
 }
\endproclaim

\enddocument